\documentclass[12pt,leqno,a4paper]{report}

\usepackage{amsmath}

\usepackage{amssymb}

\usepackage{cite}

\usepackage{dsfont}

\usepackage{fancyheadings}

\usepackage{rotating}

\usepackage{theorem}

\usepackage{verbatim}

\usepackage[all]{xy}
\usepackage{lscape}
\usepackage{currvita}

\newtheorem{definition}{Definition}
\newtheorem{lemma}[definition]{Lemma}
\newtheorem{proposition}[definition]{Proposition}
\newtheorem{theorem}[definition]{Theorem}
\newcommand{\Ort}{\ensuremath{\mathcal{I}}}
\newcommand{\Top}{\ensuremath{\mathcal{U}}}

\newcommand{\re}{\ensuremath{\mathbb{R}}}
\newcommand{\cx}{\ensuremath{\mathbb{C}}}

\newcommand{\ewedge}{\ensuremath{\overline{\wedge}}}
\newcommand{\pp}{\ensuremath{\underline{p}}}
\newcommand{\spec}{\ensuremath{\mathcal{S}}}
\newcommand{\specb}{\ensuremath{\mathcal{S}_B}}

\newcommand{\KB}{\ensuremath{\mathcal{K}_B}}

\newcommand{\Pj}{\ensuremath{\mathcal{P}(j)}}

\newcommand{\oneb}{\ensuremath{\mathbf{1}_B}}
\newcommand{\onebpu}{\ensuremath{\mathbf{1}_{BPU}}}
\newcommand{\FM}{\ensuremath{\mathbf{F}_M}}
\newcommand{\spinc}{\ensuremath{Spin^{c}}}
\newcommand{\RM}{\ensuremath{\mathbf{R}_M}}
\newcommand{\rE}{\tilde{E}^{\ast}}
\newcommand{\norm}{\ensuremath{\|.\|}}
\author{Robert Waldm\"uller 
\thanks{ Partially supported by the DFG}}

\title{Products and push-forwards in parametrised
cohomology theories} 
\begin{document}

\maketitle
\thispagestyle{empty}
\tableofcontents
\chapter*{Introduction}
\markboth{Introduction}{Introduction}

It is well known that 
generalised cohomology theories defined on topological spaces correspond 
to spectra. Under this correspondence, the cohomology groups $E^{\ast}(X)$ of a space $X$ arise 
as maps in the stable homotopy category from the suspension spectrum 
$\Sigma^\infty X $ of $X$ to $E$,
the spectrum representing the cohomology theory.
If one assumes that $X$ is a CW-complex and $E$ is an $\Omega$-spectrum
(that is, the maps $ E(n) \rightarrow \Omega_m E(n+m)$ are weak homotopy
equivalences), 
elements of the  $n$th group can be realised as maps from $X$ to
$E(n)$. 
Maps from $X$ to $E(n)$ are the same as sections
in the trivial bundle over $X$ with fibre $E(n)$.
The idea to define twisted cohomology groups
is to replace the trivial $E(n)$-bundle over $X$ by a a non-trivial $E(n)$-bundle
$M$ over $X$
and define the twisted groups as homotopy classes of sections 
of
$M$. 
Assuming that  
$\mathcal{E}(n) \rightarrow B$ is the universal  $E(n)$-bundle, there is a map
$q:X \rightarrow B$ such that $M\cong q^\ast \mathcal{E}(n)$. Sections
of $M$ are nothing but maps $u: X \rightarrow \mathcal{E}(n)$ 
such that
$$\xymatrix{X \ar[rr]^u \ar[dr]^q & & \mathcal{E}(n) \ar[dl] \\
& B &}$$ 
commutes, i.e. maps from $X$ to $M$ in the category of spaces over $B$.
This leads to the definition of a twisted or parametrised cohomology theory as a cohomology
theory on the category of spaces over a fixed base space $B$. 
In this picture of twisted cohomology theories, the set of twists on
a given space $X$ corresponds to the set of maps 
$X \rightarrow B$; we will see later
that a homotopy of twists defines an isomorphism
of the corresponding twisted cohomology groups. \\
Since a cohomology theory should have suspension isomorphisms
relating the groups in different degrees, the spaces
$\mathcal{E}(m)$ should be interrelated. 
More precisely, there should be maps
$\Sigma^m \mathcal{E}(n) \rightarrow \mathcal{E}(n+m)$, where $\Sigma$ denotes
fibrewise suspension, i.e., $\mathcal{E}$ should be what is called a parametrised
spectrum or ex-spectrum. Ex-spectra were
first studied in \cite{Clapp, ClappPuppe}. The category of
ex-spectra $\specb$ over a fixed base $B$ is similar to the  
category of non-parametrised spectra $\spec$. In particular, one can 
define the stable
parametrised homotopy category $Ho\specb$ in which 
the suspension functor is invertible. 
Moreover,  
just as  generalised cohomology theories on spaces correspond to spectra, 
cohomology theories on 
the category of spaces over a 
fixed base correspond to
ex-spectra \cite{Schoen}.  \\
In \cite{MS}, the authors discuss $Ho\specb$  
in the language of model categories. Using the concept of orthogonal spectra, they are able to define a symmetric monoidal structure which exists already before passage to the stable homotopy category. \\

An important point in the application of twisted cohomology theories 
is the existence of a product. In the examples, there is 
frequently the notion of the sum $\alpha + \beta$ of twists $\alpha,\beta$. 
The products in  the examples of multiplicative twisted cohomology theories and the sum of twists are related in that the
product adds the twists,
i.e. if $u$ is an $\alpha$-twisted class and $v$ is a $\beta$-twisted class,
$u \cup v$ is usually $\alpha + \beta$ twisted. The fact that 
multiplicative non-parametrised cohomology theories correspond to ring spectra
suggests that multiplicative parametrised cohomology theories 
should correspond to parametrised ring spectra. However,
the symmetric monoidal structure on $Ho\specb$ constructed in 
\cite{Clapp,ClappPuppe,MS} is not suitable for this purpose, 
because the product would not add twists.
In addition, if one thinks of twists as maps to some space $B$,
it is for a general $B$ not clear what the sum of 
twists should be. 
If $B$ is an $H$-space, one can use 
the product on $B$ to define the sum $f+g$ of maps $f,g: X \rightarrow B$.
To obtain our symmmetric monoidal structure on $Ho\specb$,
an $H$-space structure on $B$ is not enough. We need better homotopical
control of the product on $B$. In particular, the product should be not
only homotopy associative, but the associativity homotopy should be unique
up to higher homotopies. This is exactly what is captured by the notion of
an action of an $E_\infty$-operad on $B$. Such an action defines
products on $B$ in which all commutativity and associativity 
diagrams commute up to homotopy, and these homotopies are all unique up 
to higher homotopies. We define symmetric monoidal structures $\pp$
on $Ho\specb$ for base spaces $B$ with an action of an $E_\infty$-operad.
In fact, we only need 
up to four-fold associativities. Therefore, a truncated action would be sufficient.
 We use the symmetric monoidal strucures $\pp$
to define multiplicative
parametrised cohomology theories on the category of spaces over such a $B$. \\
\\
In the past few years, a twisted cohomology theory of particular interest was 
$H^3$-twisted $K$-theory \cite{AS1, AS2, FHT1, FHT2, FHT3, TXL1,TX}.
There are various pictures of twisted $K$-theory, 
corresponding to different pictures of 
$H^3$.\\
One way to interpret elements of $H^3(X,\mathbb{Z})$ is as principal 
$PU(\mathbf{H})$-bundles over $X$, where $\mathbf{H}$ denotes some fixed separable Hilbert space. 
In \cite{AS1, AS2}, the authors associate to a 
$PU(\mathbf{H})$-bundle a $Fred(\mathbf{H})$-bundle, using the conjugation action of $PU(\mathbf{H})$ 
on $Fred(\mathbf{H})$, and define twisted $K$-theory in terms of homotopy classes of sections
of the $Fred(\mathbf{H})$-bundle.   
Another interpretation of $H^3$ (which readily carries over to the equivariant setting)
is used in \cite{TXL1, TX}. An $H^3$-class of a topological stack corresponds
to a central $S^1$-extension of a groupoid representing the stack, and twisted
$K$-theory is defined as the $K$-theory of the reduced $C^\ast$-algebra of this extension. \\ 

Recall that a a class in $H^3(X,\mathbb{Z})$ can be realised as a
continuous map $X \rightarrow K(\mathbb{Z},3)$. Therefore, the category of 
"spaces with $H^3$-twists" may be identified with 
the category of spaces over $K(\mathbb{Z},3)$. In consequence, 
to define $H^3$ twisted $K$-theory as a multiplicative parametrised cohomology
theory amounts to constructing
a parametrised ring spetrum 
over $BPU(\mathbf{H})\cong K(\mathbb{Z},3)$ 
representing twisted $K$-theory, and we do precisely that. 
The properties of twisted $K$-theory such 
as homotopy invariance, Mayer-Vietoris sequences, excision and so forth follow
from the general statement that it is a multiplicative cohomology theory on spaces over
$BPU(\mathbf{H})$. In particular, the Thom isomorphisms and push-forwards in 
twisted $K$-theory which
were geometrically constructed in \cite{CaWang} follow from
the corresponding construction for multiplicative twisted cohomology theories.  
\\    

The text is organised as follows.
In the first chapter, we summarise the construction of $Ho\specb$ and 
some of its properties
from \cite{MS}. We define the symmetric monoidal structure $\pp$ in Proposition
\ref{symmstruct}. If a group $G$ acts on a spectrum $F$, an operad acts
on both $G$ and $F$ and the actions are compatible, we construct a 
parametrised ring spectrum with fibre $F$ over $BG$  in Lemma \ref{Monoid}.
\\
In the second chapter, we define twisted cohomology theories
and products in twisted cohomology theories. 
We prove a Thom isomorphism
theorem for generalised twisted cohomology theories and
define push-forward homomorphisms for proper oriented maps. 
In addition, we prove functoriality of push-forwards and a projection formula.
\\
In the last chapter, we construct parametrised ring spectra representing
twisted $K$-theory and twisted 
$\spinc$-cobordism. \\
Throughout the text, we use Quillens notion of a model category as defined in
\cite{Quillen}. A concise introduction to model categories is \cite{DS}, a 
detailed account with an eye to examples is \cite{Hovey}. All model categories
in our text are cofibrantly generated and therefore admit 
cofibrant and fibrant replacement functors which we denote by $Q$ and $R$.
For a topologically enriched category $\mathcal{C}$ and 
$X,Y \in Ob(\mathcal{C})$, we denote both the morphism set and the 
morphism space by $\mathcal{C}(X,Y)$.

\chapter{Basic definitions and constructions}
\pagestyle{fancyplain}
\section{Ex-spaces}
\begin{definition}
An ex-space over a topological space $B$ is
a space $X$ with continuous maps $q: X \rightarrow B$ and
$s: B \rightarrow X$ such that $q \circ s = id_B$.   
\end{definition}
The point-set topology to set up  convenient categories of ex-spaces
is quite tedious. We refrain ourselves from a lengthy discussion, but
rather refer the interested reader to \cite{Booth1, Booth2, Booth3, BrBooth1,
BrBooth2, James1,MS} and give a short account of the results we need.
Recall that a topological space $X$  is called 
a $k$-space if
a subset $A \subset X$ is closed  if $f^{-1}(A)$ is closed 
for all compact Hausdorff spaces $K$ and all continuous maps $f: K \rightarrow X$.
$X$ is called weak Hausdorff if $f(K)$ is closed for all compact Hausdorff 
spaces $K$ and all continuous maps $f: K \rightarrow X$.
Let $\Top$  denote the category of weak Hausdorff $k$-spaces and $\mathcal{K}$ the category of $k$-spaces.
\begin{definition}
For $B \in \Top$, define $\KB$ to be the category of ex-spaces $X$ over $B$
with $X \in Ob(\mathcal{K})$. Define
$\mathcal{K}/B$ to be the category of $k$-spaces over $B$.
The functor $_+: \mathcal{K}/B \rightarrow \KB$ is given by adjoing a 
disjoint basepoint, i.e. a disjoint copy of $B$ to spaces $X$ over $B$.
\end{definition}
We topologize the morphism sets with the subspace topology,
where we use the compact-open topology for $\mathcal{K}$. 
To improve legibility,
we will omit the maps $q: X \rightarrow B$ and $s: B \rightarrow X$ from our notation most of the time and only include them when necessary. 
For $K \in \mathcal{K}$ and $(X,q) \in \mathcal{K}/B$ we define
$$Map_B(K,X) := \{f: K \rightarrow X| \ \exists \ b \in B \  s.t. \ f(K) 
\subset q^{-1}(b)\} \subset Map(K,X),$$ again using the compact-open topology on 
$Map(K,X)$ and the obvious map to $B$. Together with 
$X \times_B K := (X \times K, q \circ pr_X )$ this defines a tensor and
cotensor structur of the topological
category $\mathcal{K}/B$ over $\mathcal{K}$, 
i.e. there are natural homeomorphisms
$$\mathcal{K} / B (X \times_B K, Y) \cong
\mathcal{K}(K, \mathcal{K} / B (X,Y))
\cong \mathcal{K} / B (X, Map_B (K, Y))
$$
Similarly, the pointed category $\KB$ is tensored and cotensored as topological category over 
$\mathcal{K}_{\ast}$. 
\begin{definition}
For ex-spaces $X$ and $Y$ over $B$, define the smash product\\ $X \wedge_B Y$ 
by the  following pushout in $\mathcal{K} /B.$ \\
$$
\xymatrix{
X \vee_B Y \ar[r] \ar[d] & X \times_B Y \ar[d] \\
\ast_B \ar[r] & X \wedge_B Y.
}$$
\end{definition}
There are internal mapping spaces $Map_B(.,.)$ in $\mathcal{K}/B$
and $F_B(.,.)$ in \KB, giving $\mathcal{K}/B$ the structure of a 
closed cartesian category and $\KB$ the structure of a closed symmetric
monoidal category \cite{Booth1,Booth2,Booth3}. \\
A continuous map $f: A \rightarrow B$ defines functors
$$f_{!}: \mathcal{K}_A \stackrel{\textstyle{\longrightarrow}}{\longleftarrow} 
\mathcal{K}_B : f^{\ast}$$
by setting
$$\xymatrix{
A \ar[d]^s \ar[r]^{f} & B \ar[d]^t \\
X \ar[r] &  f_! X} \qquad \qquad \qquad
\xymatrix{f^\ast Y \ar[d]^{r'} \ar[r] & Y \ar[d]^{q'} \\
A \ar[r]^{f} & B}$$ 
for $(X,q,s) \in \mathcal{K}_A$, $(Y,q',s') \in \mathcal{K}_B$.
The first diagram is a pushout and we define the missing map
$r: f_!X \rightarrow B$ to be the pushout of $id_B: B \rightarrow B$
and $f\circ q:X \rightarrow B$ ;
the second diagram is a pullback and 
$t': A \rightarrow f^\ast Y$ is defined to be the pullback of
$id_A: A \rightarrow A$ and $s' \circ f: A \rightarrow Y$. These functors are adjoint as functors
of topologically enriched categories, i.e. they are continuous and
the adjunction isomorphisms 
$$\KB(f_!X,Y) \cong \mathcal{K}_A(X,f^\ast Y)$$
are homeomorphisms  \cite[2.1.2]{MS}.

Furthermore, there are external products
$\overline{\times}: \mathcal{K}/A \ \times \ \mathcal{K}/B 
\rightarrow 
\mathcal{K}/(A \times B)$ and \\
$\overline{\wedge}:\mathcal{K}_{A} \times 
\KB \rightarrow \mathcal{K}_{A \times B}$ and mapping spaces
$\overline{Map}$ and $\overline{F}$ with natural isomorphisms 
$$ \mathcal{K}/ {(A \times B)}  (X \overline{\times } Y, Z) \cong 
\mathcal{K}/ {A}(X,\overline{Map}(Y,Z))  
$$and$$
\mathcal{K}_{A \times B}(X \overline{\wedge } Y, Z) \cong 
\mathcal{K}_{A}(X,\overline{F}(Y,Z)).$$ 

The base change functors can be used to derive the 
external products from the internal ones; for example, 
$X \overline{\wedge} Y := \pi_{A}^{\ast} 
X \wedge_{A \times B} \pi^{\ast}_{B} Y$. 
A pointed space $K \in \mathcal{K}_{\ast}$ is just an ex-space over
a point. Using the projection $\pi: B \rightarrow \ast$, the trivial ex-space with fibre $K$ is  
$$K_B := \pi^{\ast} K = K \times B \in \KB.$$ 
Notice that there are natural isomorphisms
$$X \overline{\wedge} K \cong X \wedge_B K \cong X\wedge K_B,$$
relating the tensor structure of $\KB$ over $\mathcal{K}$, the
external and the internal smash product; 
see \cite{MS} for a more thorough discussion
of the functors \\ $\overline{\wedge}, \wedge,
f^\ast, f_!, \wedge_B, Map_B, F_B, \overline{F}$ and their relations.      
\begin{proposition}\cite[2.2.1]{MS}
For a map $f: A \rightarrow B$, the functor 
$f^{\ast}: \KB \rightarrow \mathcal{K}_A$ is closed symmetric monoidal.
\end{proposition}

In \cite{MS}, the authors develop model structures on $\mathcal{K}/B$ and 
$\KB$. We will now give a short summary of their results, again omitting the proofs.

\begin{definition}
$f \in \mathcal{K}/B (X,Y)$ is a fibrewise cofibration
 if it has the left lifting property with respect to the maps
$$\pi_0: Map_B(I,Z) \rightarrow Z  \ \forall Z \in \mathcal{K}/B,$$
where $\pi_0$ is given by evaluation at zero.
\end{definition}
Denote the upper (lower) hemisphere of a $n$-disc by $S^{n-1}_{+} \ (S^{n-1}_{-})$
with the inclusion $i_{\pm}: S^{n-1}_{\pm} \rightarrow D^{n}$.
Similarly, we write
$$i: S^{n-1} \rightarrow D^n \qquad \text{and} \qquad j_{-}:S^{n-1} \rightarrow
S^{n}_{-}$$ for the inclusions. 
\begin{definition}
\begin{eqnarray*}
I_B & := &\{i: (S^{n-1},d \circ i) \rightarrow (D^n,d)\ |\ i\ is\ a\ fibrewise\ 
cofibration\}\\
J_B & := & \{i_{+}: (S^{n}_+,d \circ i_{+}) \rightarrow (D^{n+1},d)\ |\ 
i: (S^{n},d \circ i) \rightarrow (D^{n+1},d) \
 and  \\
& & j_{-}: (S^{n-1},d \circ i_{-} \circ j_{-}) \rightarrow (S^{n}_{-}, d \circ i_{-})\
are \  fibrewise \ 
cofibrations \}
\end{eqnarray*}
\end{definition}

\begin{theorem}\cite[6.2.5]{MS}
$\mathcal{K} / B$ is a cofibrantly generated model category,
with weak equivalences given by weak homotopy equivalences of 
total spaces and $I_B$ $(J_B)$ as set of generating (acyclic) cofibrations.
\end{theorem}

\begin{theorem}\cite[6.2.6]{MS}\label{modelkb}
$\KB$ is a cofibrantly generated model category,
with weak equivalences given by weak homotopy equivalences of 
total spaces and $(I_B)_+$ $((J_B)_+)$ as set of generating (acyclic) cofibrations.
\end{theorem}
\begin{proposition}\cite[7.3.4, 7.3.5]{MS}
For a continuous map $f: A \rightarrow B$, the adjunction $(f_!,f^*)$ is a Quillen adjunction between $\mathcal{K}_A$ and $\KB$.  
It is a Quillen equivalence if $f$ is a weak homotopy equivalence.
\end{proposition}
The internal smash product $\wedge: \KB \times \KB \rightarrow \KB$ is
not a Quillen bifunctor. Nonetheless, the authors of \cite{MS} construct a
"derived" internal smash product on the homotopy category.
However,
the external smash product 
$$\overline{\wedge}: \mathcal{K}_A \times \KB \rightarrow 
\mathcal{K}_{A \times B}$$ is a Quillen bifunctor \cite[7.3]{MS} and
thus can be derived to obtain
$$L \overline{\wedge}: Ho \mathcal{K}_A \times Ho \KB \rightarrow 
Ho \mathcal{K}_{A \times B}.$$  

For an ex-space $E \in \mathcal{K}_A$, there is a natural map
$$i_{0 !}QE \cong QE \coprod_{A \times \{0\}} A \times I
\tilde{\rightarrow} RE \times I \cong \pi_{A}^\ast RE  
$$
given by the extension of $Q E \rightarrow R Q E \rightarrow R E$. 
This is a weak equivalence,
hence it induces a natural isomorphism of total derived functors
$L i_{0 !} \cong R \pi_{A}^{\ast}$. Similarly, 
$L i_{1 !} \cong R \pi_{A}^{\ast}$ and therefore we get a natural isomorphism
$\rho : L i_{0 !} \tilde{\rightarrow} L i_{1 !}$.

\begin{proposition}\label{Htpyisokb}
A homotopy $h: A \times I \rightarrow B$
between  $f,g : A \rightarrow B$ gives rise to a natural isomorphism
$\phi: Lf_{!} \ \tilde{\rightarrow} \ Lg_{!}$.\\
If  two homotopies $h_0, h_1$ are homotopic (as homotopies)
 via $H$, the isomorphims $\phi_0$ and
$\phi_1$ coincide.
\end{proposition}

\textbf{Proof:}
The natural isomorphisms 
$$h_!  i_{0!} \cong f_!, \qquad h_!  i_{1!} \cong g_!$$
between left Quillen functors yield isomorphisms between the derived functors 
\cite{Hovey}. We combine these and the natural isomorphism
$\rho:L i_{0 !} \tilde{\rightarrow} L i_{1 !}$ to obtain the desired 
isomorphism  
$$Lf_{!} \cong Lh_! \ L i_{0!} \stackrel{Lh_{!}(\rho)}{\longrightarrow} Lh_{!} \ L i_{1 !} \cong L g_!.$$ 
Let us now discuss the case in which there are two
homotopic homotopies $h_0, h_1$ from $f$ to $g$.  
For the inclusions $j_0, j_1: A \times I \rightarrow A \times I \times I$,
and $E \in \mathcal{K}_A$,
we have again  natural weak equivalences in $\mathcal{K}_{A \times I \times I}$
$$  j_{0 !} Qi_{0 !} QE  \stackrel{\alpha_0}{\rightarrow} \pi_{A \times I}^\ast
R  i_{0 !} QE \stackrel{\alpha_1}{\leftarrow}j_{1 !} Qi_{0 !} QE$$
$$  j_{0 !} Qi_{1 !} QE  \stackrel{\beta_0}{\rightarrow} \pi_{A \times I}^\ast
R  i_{1 !} QE \stackrel{\beta_1}{\leftarrow}j_{1 !} Qi_{1 !} QE$$
$$ j_{0 !} Q\pi_{A}^{\ast}R E  \stackrel{\gamma_0}{\rightarrow} 
\pi_{A \times I}^\ast
R Q \pi_{A}^{\ast} R E \stackrel{\gamma_1}{\leftarrow}j_{1 !} 
Q\pi_{A}^{\ast}RE.$$
Note that $\pi_{A \times I}^{\ast} \cong \overline{\wedge}S^{0}_{I}$, 
hence it preserves cofibrant objects.
For $E \in \mathcal{K}_A$, commutativity of the diagram
$$
\xymatrix{
 &  h_{0 !} Q  i_{0 !} Q E \ar[r]^{L h_{0 !} (\rho)} \ar[d]^{H_! (\alpha_0)} 
& h_{0 !} Q i_{1 !} Q E
\ar[dr]^{~} \ar[d]^{H_! (\beta_0)} \\
 f_! Q E \ar[ur]^{~}  \ar[dr]^{~}& 
H_! \pi_{A \times I}^\ast R  i_{0 !} Q E \ar[d]^{H_! (\alpha_1)^{-1}}&
H_! \pi_{A \times I}^\ast R  i_{1 !} QE \ar[d]^{H_! (\beta_1)^{-1}} &  
g_! QE \\
&  h_{1 !} Q  i_{0 !} Q  E \ar[r]^{L h_{1 !} (\rho)}  
&  h_{1 !}  Q i_{1 !} Q E \ar[ur] 
}$$
in $Ho\KB$ gives the result, so 
let us first show that the middle square commutes. 
A lift of the square  to $\KB$ is $H_!$ applied to the following diagram in
$\mathcal{K}_{A \times I \times I}.$ 
$$\xymatrix{
j_{0 !}Qi_{0!}QE \ar[r]\ar[d]^{\alpha_0} & 
j_{0 !}Q \pi_{A}^{\ast} RE \ar[d]^{\gamma_0}& 
j_{0 !}Qi_{1!}QE \ar[l] \ar[d]^{\beta_0}\\
 \pi_{A \times I}^{\ast} R  i_{0!} QE \ar[r] & 
 \pi_{A \times I}^{\ast} R Q \pi_{A}^{\ast} RE &
 \pi_{A \times I}^{\ast} R  i_{1!} QE \ar[l] \\
j_{1 !}Qi_{0!}QE \ar[r] \ar[u]^{\alpha_1} & 
j_{1 !}Q \pi_{A}^{\ast} RE \ar[u]^{\gamma_1}& 
j_{1 !}Qi_{1!}QE. \ar[l] \ar[u]^{\beta_1}
}$$ 
The top and bottom horizontal maps stem from the natural maps
$$i_{0 _!} QE \rightarrow \pi_{A}^{\ast} RE \qquad  \qquad
i_{1 _!} QE \rightarrow \pi_{A}^{\ast} RE.$$ Since 
$i_{0 !} QE, \ i_{1!}QE$ are cofibrant, these maps can be lifted
to maps $$i_{0 _!} QE \rightarrow Q \pi_{A}^{\ast} RE \qquad  \qquad
i_{1 _!} QE \rightarrow Q \pi_{A}^{\ast} RE$$ which yield the middle
horizontal arrows.
The above diagram commutes in $\mathcal{K}_{A \times I \times I}$ and all maps are weak equivalences. Since  all objects in the diagram are cofibrant and 
$H_!$, being a left Quillen functor, preserves all weak equivalences 
between cofibrant objects, the middle square commutes in $Ho \KB$.
  
To show commutativity of the left triangle, we choose again a lift to 
$\mathcal{K}_B$, namely
$$\xymatrix{
 & h_{0!} Q i_{0!} Q E \ar[d]^{H_!(\alpha_0)} \ar[dl]_{h_{0!} (q)}  \\
f_! QE & H_{!}\pi_{A \times I}^{\ast} R  i_{0!} QE \\
 & h_{1!} Q i_{0!} Q E \ar[u]_{H_!(\alpha_1)} \ar[ul]^{h_{1!} (q)},  
}$$
with $q$  being the natural map $ Q i_{0!} Q E \rightarrow i_{0!} Q E$.
Observe that we have factorisations
$$\alpha_0 : j_{0 !} Qi_{0 !} QE \rightarrow \pi_{A \times I}^\ast i_{0 !} QE \rightarrow
\pi_{A \times I}^\ast
R  i_{0 !} QE $$
$$\alpha_1: j_{1 !} Qi_{0 !} QE
\rightarrow \pi_{A \times I}^\ast i_{0 !} QE \rightarrow
\pi_{A \times I}^\ast
R  i_{0 !} QE. $$
We use these to see that
$H_!(\alpha_0)$ and $H_!(\alpha_1)$ factor through
$H_! \pi_{A \times I}^{\ast}i_{0!} QE$ and we get a commutative diagram
$$\xymatrix{
 & h_{0!} Q i_{0!} Q E \ar[dr]^{H_!(\alpha_0)} \ar[dl]_{h_{0!} (q)} \ar[d] &   \\
f_! QE & H_! \pi_{A \times I}^{\ast}i_{0!} QE \ar[r] & 
H_{!}\pi_{A \times I}^{\ast} R i_{0!} QE \\
 & h_{1!} Q i_{0!} Q E \ar[ur]_{H_!(\alpha_1)} \ar[ul]^{h_{1!} (q)} \ar[u]&  
}.$$

Furthermore, there are maps
$$\xi_s: j_{s!} i_{0!} QE \rightarrow \pi_{A \times I}^{\ast}i_{0!} QE \ \ \ 
\forall \, s \, \in I
$$
given by extension of the identitiy $i_{0!} QE \rightarrow i_{0!} QE$.
The triangles
$$\xymatrix{
&  & h_{0!} Q i_{0!} Q E \ar[dll]_{h_{0!} (q)} \ar[d]    \\ 
f_! QE \ar[rr]^{H_! (\xi_0)} &   & H_! \pi_{A \times I}^{\ast}i_{0!} QE 
}
$$
and 
$$\xymatrix{
f_! QE \ar[rr]^{H_! (\xi_1)} &   & H_! \pi_{A \times I}^{\ast}i_{0!} QE \\
& & h_{1!} Q i_{0!} Q E \ar[ull]^{h_{1!} (q)} \ar[u]
}
$$
commute. Since $H \circ j_s \circ i_0 = f \,\, \, \forall \, s  \in  I$, 
$H_!(\xi_0)$ and $H_!(\xi_1)$ are homotopic over $B$. Thus,
the left triangle commutes in $Ho\KB$.\\

The argument for the right triangle is similar.
\hfill $\square$\\ \\
In addition, the correspondence between homotopies
and natural equivalences is functorial in the following sense.
\begin{proposition}\label{funckb}
Let $h_0$ be a homotopy between $f_0, f_1: A \rightarrow B$ and $h_1$ be one
connecting $f_1$ and $f_2$, and denote by $h$ the homotopy from 
$f_0$ to $f_2$ given by composing $h_0$ and $h_1$.
 Then the corresponding isomorphisms
$\phi_0: Lf_{0!} \tilde{\rightarrow} Lf_{1!}$, 
$\phi_1: Lf_{1!} \tilde{\rightarrow} Lf_{2!}$
and $\phi: Lf_{0!} \tilde{\rightarrow} Lf_{2!}$
satisfy
$\phi_1 \circ \phi_0 = \phi: Lf_{0!}  \tilde{\rightarrow} Lf_{2!}$.
\end{proposition}
\textbf{Proof:}
Consider the double interval $ J:=[0,2]$ and the maps
$$ 
\pi: A \times J \rightarrow A, \  
j_0: A \times I \rightarrow A \times J , \ 
j_1 : A \times I \rightarrow A \times J, \ i_{kl}:=j_{k} \circ i_{l}$$
and $h$ as a map $A \times J \rightarrow B$. 
There are identities
$$ f_{k+l}= h_k \circ i_l = h \circ i_{kl} \ \forall k,l \in \{0,1\}.$$
Consider the diagram
$$\xymatrix{
& Q j_{0!} i_{0 !} Q E \ar[d] \ar[dr] & \\
j_{0!} Q i_{0 !} QE \ar[r] \ar[d] \ar[ur]^{\alpha_2}&
 j_{0!} i_{0!} QE \ar[d]  
\ar[r]^{\alpha_0} & 
Q\pi^\ast RE \ar[d] \\
j_{0!} Q \pi_{A}^{\ast} RE \ar[urr]^{\alpha_1} \ar[r] &
j_{0!}  \pi_{A}^{\ast} RE \ar[r] & \pi^\ast RE}
$$
in which all maps except $\alpha_i$ are either natural
"extensions by zero" or of the form $QF \rightarrow F$, and
$\alpha_i$ are obtained by lifting the respective maps to
the cofibrant replacements. 
Since maps of the form $QF \rightarrow F$ are weak equivalences,
the diagram commutes in $Ho \mathcal{K}_{A \times J}$.
Similarly, let $\alpha_3$ be defined to be the lift in
$$\xymatrix{
 & 
 &
Q \pi^\ast RE \ar[d]\\
j_{0!} i_{1!} QE \ar[r] \ar[rru]^{\alpha_3}& 
j_{0!} \pi_{A}^{\ast} RE \ar[r] &
 \pi^\ast RE }$$
and observe that 
$$\xymatrix{
j_{0!} Q i_{1!} QE \ar[r] \ar[d]& 
j_{0!} Q \pi_{A}^{\ast} RE \ar[r]^{\alpha_1}  &
Q \pi^\ast RE \\
j_{0!} i_{1!} QE \ar[rru]^{\alpha_3}& 
}$$
commutes in $Ho\mathcal{K}_{A\times J}$.
Putting everything together, 
applying $h_!$ and using the natural identifications
$$h_! \circ j_{k!} \circ i_{l!} \cong f_{k+l !},$$ we obtain
$$\xymatrix{
 h_!Q i_{00!}QE \ar[rrr] & & &  h_! Q\pi^{\ast} RE \\  
& h_{0!}Q \pi_{A}^{\ast} RE \ar[urr]^{h_!(\alpha_1)}  \\
 h_{0!} Q i_{0!} QE \ar[uu]^{h_!(\alpha_2)} \ar[ur] \ar[d] & & 
h_{0!} Q i_{1!} QE \ar[ul] \ar[dr] &  
 \\
f_{0!}QE  &  & & f_{1!}QE \ar[uuu]^{h_!(\alpha_3)}.
}$$
One can repeat the construction with $h_1$ instead of $h_0$ and the 
second half of $J$. This completes the previous diagram to  
$$\xymatrix{
h_!Q i_{00!}QE \ar[rr] & &  h_! Q\pi^{\ast} RE & &  
h_!Q i_{11!}QE \ar[ll]\\  
& h_{0!}Q \pi_{A}^{\ast} RE \ar[ur] &  & h_{1!}Q \pi_{A}^{\ast} RE \ar[ul] & \\
h_{0!} Q i_{0!} QE \ar[uu] \ar[ur] \ar[d] & 
h_{0!} Q i_{1!} QE \ar[u] \ar[uur] \ar[dr] &  
& h_{1!} Q i_{0!} QE \ar[uul] \ar[u] \ar[dl]&
  h_{1!} Q i_{1!} QE \ar[uu] \ar[ul] \ar[d] \\
f_{0!}QE  & & f_{1!}QE \ar[uuu] & & f_{2!} QE.
}$$
\\
The natural transformation corresponding to $h$ is obtained by the left and 
right vertical arrows and the top horizontal ones whereas the $\phi_i$ correspond
to the left
and right 'hats'.    
\hfill$\square$
\section{Ex-spectra}
\begin{definition}
Let $\mathcal{J}$ denote the category of finite dimensional
real inner product spaces and linear isometric isomorphisms.
A $\mathcal{J}$-space  over $B \in \Top$ is a continuous functor
$E: \mathcal{J} \rightarrow \KB$; $\mathcal{J}\KB$ is the category of  
$\mathcal{J}$-spaces  over $B$.
\end{definition}
The 
external smash product of $X,Y \in\mathcal{J}\KB$ is given by
$$X \overline{\wedge}_{B} Y := \wedge \circ (X \times Y): \mathcal{J} \times \mathcal{J}
\rightarrow \KB.$$
Define the smash product $X\wedge_B Y$ to be the topological left 
Kan extension of $X \overline{\wedge}_{B} Y$ along 
$\oplus:\mathcal{J}\times \mathcal{J} \rightarrow \mathcal{J}$.

\begin{theorem}\cite[11.1.3, 11.1.6]{MS} 
$(\mathcal{J}\KB,\wedge_B)$ is a 
closed symmetric monoidal category. It is tensored and cotensored over
\KB.
\end{theorem} 
The tensor and cotensor structure are defined levelwise.
Moreover, for a map $f: A \rightarrow B$ of base spaces,
levelwise application of the adjunction
$(f_!,f^\ast)$ defines an adjunction
$$f_! : \mathcal{J}\mathcal{K}_A \stackrel{\textstyle{\rightarrow}}{\leftarrow} 
\mathcal{J}\KB: f^\ast.$$
For a real inner product space $V$, let $S^V$ denote the one-point compactification.
Define the sphere $S_B$ over $B$ to be the $\mathcal{J}$-space sending
$V$ to $S^{V}_B $. 
This is a commutative monoid in $\mathcal{J}\KB$,
so we can define  $\specb$,  the ex-spectra over $B$, to be the 
$S_B$-modules. The smash product $X\wedge Y$ of ex-spectra $X,Y$ is the
coequalizer (in $\mathcal{J}\KB$) of
$$ X \wedge_B S_B \wedge_B Y 
\stackrel{\textstyle{\longrightarrow}}{\longrightarrow} X \wedge_B Y.$$  

\begin{proposition}\cite[11.2.5]{MS} 
 $(\specb,\wedge)$ is a 
closed symmetric monoidal category with unit $S_B$.
\end{proposition} 
Again, levelwise application of the functors
$f_!,f^\ast$ gives an adjunction \\
$f_!: 
\mathcal{S}_A \stackrel{\textstyle{\rightarrow}}{\leftarrow} \specb : f^\ast$
 and the functor
$f^\ast: \specb \rightarrow \mathcal{S}_A$ is closed symmetric monoidal. 
\begin{definition}\cite{MS}For base spaces $A,B$, the external smash product is
 $$\overline{\wedge} := 
\wedge \circ (\pi_{A}^{\ast} \times \pi_{B}^{\ast}): \mathcal{S}_A \times \specb \rightarrow \mathcal{S}_{A \times B}.$$ 
\end{definition}
There is an obvious pair of adjoint functors 
$$\Sigma^{\infty}: \KB \rightarrow \specb \qquad \text{ and } \qquad
 \Omega^{\infty}: \specb \rightarrow \KB,$$ given by
$\Sigma^{\infty}(X)(V)= X \wedge_B S_{B}(V)$ and 
$\Omega^{\infty}(E):=E(0)$.
More generally, the functor $\Omega_{V}^{\infty}$ given by evaluation at $V$ has a left adjoint $\Sigma_{V}^{\infty}$ \cite{MS}.
\begin{proposition}\cite[11.2.5]{MS}
$\specb$ is 
tensored and cotensored over \KB. For  $X \in \specb$ and $K \in \KB$, there are natural isomorphims
$$ X \wedge_B K \cong X \wedge_B \Sigma^{\infty} K. $$
\end{proposition}
\section{Model structures on ex-spectra}
In this section, we give a short summary 
of the model structures on ex-spectra developed in \cite{MS}. 
First, a level model structure is defined which is then utilised to
obtain a stable model structure, where the adjectives "level" and "stable" refer to
the respective classes of weak equivalences. 
\begin{definition}\cite{MS}
A map $f: X \rightarrow Y$ of spectra over $B$ is called a
level weak equivalence if $f(V): X(V) \rightarrow Y(V)$ is a 
weak homotopy equivalence for all $V \in \mathcal{J}$.
It is called a level fibration if each $f(V)$ is a fibration
in the model structure on $\KB$ defined in Theorem \ref{modelkb}. 
\end{definition}
Fix a skeleton $skel(\mathcal{J})$ of $\mathcal{J}$. 
\begin{definition}\cite{MS}
\begin{eqnarray*}
FI_B &:= &\{\Sigma_{V}^{\infty} i \, | \, i \in I_B, \ V \in  skel(\mathcal{J})\} \\
FJ_B &:= &\{\Sigma_{V}^{\infty} j  \, | \,  j \in J_B, \ V \in  skel(\mathcal{J})\} 
\end{eqnarray*}
\end{definition}
\begin{theorem}\cite[12.1.7]{MS}
$\specb$ is a cofibrantly generated topological model category with
level weak equivalences as weak equivalences,
level fibrations as fibrations and $FI_B (FJ_B)$ as generating (acyclic) cofibrations .
\end{theorem}
Recall that the homotopy groups $\pi_q (X)$ of an (unparametrised) spectrum $X$ are defined as 
the colimits of the groups $\pi_q(\Omega^V(X(V)))$.
\begin{definition}\cite{MS}
A map $f: X \rightarrow Y$ of parametrised 
spectra is called a stable weak equivalence
if, after level fibrant approximation $R^l$, 
it induces
an isomorphism $i_{b}^{\ast}R^l f: i_{b}^{\ast} R^l X \rightarrow   
i_{b}^{\ast} R^l Y$ on homotopy groups of (unparametrised) spectra
for all points $b \in B$.
\end{definition}
A map of ex-spectra is called  s-cofibration if it is a cofibration
in the level model structure. 
\begin{theorem}\cite[12.3.10]{MS}
$\specb$ is a cofibrantly generated model category
with the stable weak equivalences as weak equivalences and the s-cofibrations
as cofibrations. 
\end{theorem}
As in the unparametrised case, any level weak equivalence is a stable weak
equivalence.
The stable model structure from the above theorem is indeed stable in the sense
that the suspension functor 
$\Sigma: \specb \rightarrow \specb, \Sigma(E)= E\wedge_B  S^1 $
is invertible in the homotopy category. 
\begin{definition}\cite{MS}
An $\Omega$-spectrum over $B$ is a level fibrant spectrum over $B$ such that
each of its adjoint structure maps is a weak equivalence.
\end{definition}
It is shown in \cite[12.3.10,12.3.14]{MS} that the $\Omega$-spectra over $B$ are
the fibrant objects in the stable model structure. 
Moreover, the sets of level and stable weak equivalences between 
$\Omega$-spectra coincide.

\begin{proposition}\cite[12.6.2]{MS}
The adjunction
$\Sigma_{V}^{\infty}: \KB \stackrel{\textstyle{\rightarrow}}{\leftarrow} \specb : \Omega_{V}^{\infty}$ 
is a Quillen adjunction for all $V \in \Ort.$
\end{proposition}
\begin{proposition}\cite[12.6.7]{MS}
For a continuous map $f: A \rightarrow B$, the adjunction $(f_!,f^*)$ is a Quillen adjunction.  
It is a Quillen equivalence if $f$ is a weak homotopy equivalence.
\end{proposition}
The external smash product 
$\overline{\wedge}: \mathcal{S}_A \times \mathcal{S}_B \rightarrow 
\mathcal{S}_{A \times B} $
is homotopically well-behaved in the stable model structure and compatible
with the base change functors; more precisely, we have the following proposition.
\begin{proposition} \cite[12.6.5,13.7.2]{MS}\label{bchange} 
$\overline{\wedge}: \mathcal{S}_A \times \mathcal{S}_B \rightarrow 
\mathcal{S}_{A \times B}$
is a Quillen bifunctor.
For maps
$f:A \rightarrow B$ and $g: A' \rightarrow B'$ ,
there are natural equivalences
$$L( \overline{\wedge} \circ (f_! \times g_!))  \cong L((f \times g)_! \circ  \overline{\wedge}  ).$$   
\end{proposition}
For a spectrum $E$ over $A$, there is a natural map 
$i_{0!}Q E \rightarrow i_{0!}RQE \rightarrow \pi_{A}^{\ast} RE$ of 
spectra over $A \times I$.
The first map is a stable weak equivalence since $i_{0!}$ is a Quillen left adjoint,
and $RQE \rightarrow RE$ is a stable 
weak equivalence between fibrant spectra, hence 
it is a level weak equivalence which implies that the second map
is a (level and therefore stable) weak equivalence as well.
Moreover, since the base change functors are defined levelwise,
the same proofs as for Proposition \ref{Htpyisokb} and
\ref{funckb}
give the corresponding statements for $Ho \specb$.
\begin{proposition}\label{htpyisosb}
A homotopy $h: A \times I \rightarrow B$
between  $f,g : A \rightarrow B$ gives rise to a natural isomorphism
$\phi: Lf_{!} \ \tilde{\rightarrow} \ Lg_{!}$ of functors 
$Ho \mathcal{S}_A \rightarrow Ho \specb$.
If  two homotopies $h_0, h_1$ are homotopic (as homotopies)
 via $H$, the isomorphims $\phi_0$ and
$\phi_1$ coincide.
\end{proposition}

\begin{proposition}\label{funcsb}
Let $h_0$ be a homotopy between $f_0, f_1: A \rightarrow B$ and $h_1$ be one
connecting $f_1$ and $f_2$, and denote by $h$ the homotopy from 
$f_0$ to $f_2$ given by composing $h_0$ and $h_1$.
 Then the corresponding isomorphisms
$\phi_0: Lf_{0!} \tilde{\rightarrow} Lf_{1!}$, 
$\phi_1: Lf_{1!} \tilde{\rightarrow} Lf_{2!}$
and $\phi: Lf_{0!} \tilde{\rightarrow} Lf_{2!}$
satisfy
$\phi_1 \circ \phi_0 = \phi: Lf_{0!}  \tilde{\rightarrow} Lf_{2!}$.

\end{proposition} 
\section{Operads}
\begin{definition}
Let $\Delta$ denote the category 
with  the sets $(n)=\{0,...,n\}, n \geq 0$ as objects 
and the non-decreasing functions as morphisms.
The category of simplicial spaces  $\mathbf{S}\Top$ is 
defined as the category
of functors $\Delta^{op} \rightarrow \Top$. 
\end{definition}

Using the functor $\triangle: \Delta \rightarrow \Top$ mapping
$(n)$ to the topological n-simplex
$$
\triangle^n:=\{(t_0,...,t_n) \in \mathbb{R}^{n+1} \, | \, \Sigma t_i = 1, 
\ t_i \geq 0 \}$$
we define the geometric realisation $| \cdot | : \mathbf{S}\Top \rightarrow \Top$ 
as $$|X_\bullet|:= \int^{(n) \in \Delta} X_n \times \triangle^n ,$$
where the coend is defined as the initial object
of $X_\bullet \times \triangle^\bullet/\Top$, the category whose objects are
spaces $Y \in \Top$ together with maps $X_n \times \triangle^n \rightarrow Y$ for all $n$ such that
$$\xymatrix{
X_k \times \triangle^{l} \ar[r] \ar[d] & X_k \times \triangle^{k} \ar[d] \\
X_l \times \triangle^{l} \ar[r] & Y
}$$ commutes for all morphisms in $\Delta((l),(k))$.

\begin{proposition}\cite[11.5]{GeomILS}
Geometric realisation preserves finite products, i.e. there are natural
homeomorphisms 
$\xi: |X_{\bullet} \times Y_{\bullet}| \cong |X_{\bullet}| \times |Y_{\bullet}|$ for all
simplicial spaces $X_{\bullet}, Y_{\bullet} \in \mathbf{S}\Top$.
\end{proposition}

\begin{definition}\cite{Markl, GeomILS}
An operad (in \Top) is a sequence of spaces \\
$\mathcal{P}(j) \in \Top, \ j \in \mathbb{N}$ with right 
$\Sigma_j$-actions on $\mathcal{P}(j)$, continuous $\Sigma$-equivariant maps
$$\gamma: \mathcal{P}(k) \times \mathcal{P}(j_1) \times ... \times 
\mathcal{P}(j_k) \rightarrow \mathcal{P}(j_1 + j_2 + ... j_k)$$
satisfying associativity and a unit $1 \in \mathcal{P}(1)$.
A morphism of operads $\theta: \mathcal{P} \rightarrow \mathcal{P}'$ is a sequence of continuous $\Sigma_j$-equivariant maps 
$ \theta(j): \mathcal{P}(j) \rightarrow \mathcal{P}'(j)$ such that
$$\theta(j_1 + ... + j_k)\circ \gamma = \gamma' \circ (\theta(k) 
\times \theta(j_1) \times ... \times \theta(j_k)).$$
An $E_\infty$-operad is an operad $\mathcal{P}$ such that all
spaces $\mathcal{P}(j)$ are contractible.
\end{definition}
An operad is called pointed if it is indexed on $\mathbb{N}_0$ and  
$\mathcal{P}(0)$ consists of a single point. \\
\\
\textbf{Example:}
Let $(\mathcal{C},\wedge)$ be a 
symmetric monoidal category enriched over $\Top$. For any
 $X \in \mathcal{C}$, the endomorphism operad $End_X$ is
defined as 
$$End_X(j):= \mathcal{C}(X^{\wedge j},X) \qquad \gamma(f,g_1,...,g_k) :=
f \circ (g_1 \wedge ... \wedge g_k) $$ with the identity map as unit. 

\begin{definition}
An action of an operad $\mathcal{P}$ on $X \in \mathcal{C}$ is a morphism
of operads $\theta: \mathcal{P} \rightarrow End_X$. $(X,\theta)$ is called 
an algebra over $\mathcal{P}$. We denote the category of algebras over
$\mathcal{P}$ in $\mathcal{C}$ by $\mathcal{P}[\mathcal{C}]$.
\end{definition}

\begin{proposition}\label{symmstruct}
If a pointed $E_\infty$-operad $\mathcal{P}$ acts
on a space $B$, any $p \in \mathcal{P}(2)$
endows $Ho\specb$ with the structure of a symmetric
monoidal category.
\end{proposition}
\textbf{Proof:}
For any $q \in \mathcal{P}(j)$, 
the $\mathcal{P}$-algebra structure yields a
map $\theta(q): B^j \rightarrow B$
and  any path in $\mathcal{P}(j)$ from $q$ to $q'$ yields
a homotopy from $\theta(q)$ to $\theta(q')$.
By Proposition \ref{htpyisosb}, this homotopy 
yields a natural transformation
$$L \theta (q)_! \tilde{\rightarrow} L \theta (q')_!$$
and
since $\mathcal{P}(j)$ is contractible by assumption,
any two such paths are homotopic and thus define the same natural transformation.\\
We define the bifunctor  associated to $p \in \mathcal{P}(2)$ 
to be
$$\underline{p}:=L (\theta(p)_! \circ \overline{\wedge}): 
Ho\specb \times Ho\specb \rightarrow Ho\specb.$$
and the unit by $\oneb := i_{\ast !}S$ with 
$$i_{\ast}:=\theta(\mathcal{P}(0)): \ast \rightarrow B$$ 
and $S$ the sphere spectrum. Observe that $S$ and hence also $\oneb$ is cofibrant.
The idea to obtain unit, symmetry and associativity isomorphisms 
for $\pp$
is very simple: the things one wants to relate are obtained as push-forwards
along maps corresponding to different points in the operad.
These points in $\mathcal{P}(j)$ can be connected by paths in $\mathcal{P}(j)$
and we use the natural transformations corresponding to these paths to construct the
isomorphisms.  \\
Throughout the proof, we suppress the associativity, unit and commutativity
isomorphisms for $\overline{\wedge}$ from the notation.\\
To get the (left) unit isomorphism, observe that the map
$$p_e:  B \rightarrow B, \ \ \   b \mapsto \theta(p) (\ast,b)$$ is given by
$\theta(\gamma(p;\ast,id))$. Since $\mathcal{P}(1)$ is connected, $p_e$ 
is homotopic
to the identity\cite{Mayarticle}. 
Combining Propositions \ref{bchange} and \ref{htpyisosb},
we obtain $$\pp(\oneb,X) 
\cong \theta(p)_!( (i_{\ast !} S) \overline{\wedge} QX) \cong
\theta(p)_!(i_\ast \times id)_! S \overline{\wedge} QX \cong
p_{e!} QX \cong X.$$ Similarly, we make use of
the fact that the maps $p_{12}:=\theta(\gamma(p;p,1)): B^3 \rightarrow B$ and 
$p_{21}:=\theta(\gamma(p;1,p)): B^3 \rightarrow B$
are homotopic. \\
\begin{eqnarray*}
\pp(\pp(X,Y),Z) & \cong & \theta(p)_! (\theta(p)_! 
(QX \overline{\wedge} QY)\overline{\wedge} QZ)\\
&\cong& \theta(p)_! (\theta(p) \times id)_! 
((QX \overline{\wedge} QY) \overline{\wedge} QZ) \\
& \cong & p_{12!}( (QX\overline{\wedge}QY)\overline{\wedge}QZ ) \cong
p_{21!}( QX\overline{\wedge}(QY\overline{\wedge}QZ) ) \\
&\cong& \pp(X,\pp(Y,Z))
\end{eqnarray*} 
which yields the associativity transformation.
For the non-trivial element $\tau \in \Sigma_2$, the maps 
$\theta(p)$ and $\theta(p \tau)$ are homotopic, 
hence Proposition \ref{htpyisosb}
gives  the commutativity isomorphism $T$.\\
What remains to be checked are the usual coherence diagrams which can be
found for example in \cite{Kelly}.
 
For spectra $X,Y,W,Z \in \specb$ we need to show
that
$$\xymatrix{
\pp(X,\pp(Y,\pp(Z,W))) \ar[d] \ar[r] &
\pp(\pp(X,Y),\pp(Z,W)) \ar[r] & \pp(\pp(\pp(X,Y),Z),W)  \\
\pp(X,\pp(\pp(Y,Z),W)) \ar[rr] & & 
\pp(\pp(X,\pp(Y,Z)),W) \ar[u]
}$$
commutes. Using Propositions \ref{bchange} and \ref{htpyisosb} and writing
$$E:=QX\overline{\wedge} 
QY \overline{\wedge} QZ \overline{\wedge} QW,$$
this reduces to the commutativity of
$$\xymatrix{
p_{1!} \ E \ar[r] \ar[d]&
p_{2!} \ E
 \ar[r] & 
p_{3!}  \ E
\\
p_{4!} \ E \ar[rr] & & 
p_{5!} \ 
E \ar[u]
}$$
where the $p_i$ denote the maps corresponding
to points in $\mathcal{P}(4)$ 
(for example, $p_1 = \theta(\gamma(p;1,\gamma(p;1,p))$)
and the arrows correspond to natural isomorphisms induced
from homotopies (i.e., paths in $\mathcal{P}(4)$)
between the maps. 
$\mathcal{P}(4)$ is contractible, so all paths connecting two points
are homotopic. This observation, together with Proposition \ref{funcsb},
gives the result.
By the same argument,
$$\xymatrix{
\pp(X,Y)\ar[r] \ar[dr]^{id} & \pp(Y,X) \ar[d]\\
& \pp(X,Y) } \hspace{2cm}  
\xymatrix{
\pp(\oneb,\pp(X,Y)) \ar[r] \ar[d] & \pp(X,Y) \\
\pp(\pp(\oneb,X),Y) \ar[ur]}
$$
and 
$$\xymatrix{
\pp(X,\pp(Y,Z)) \ar[r] \ar[d] & \pp(\pp(X,Y),Z) \ar[r] & \pp(Z,\pp(X,Y)) \ar[d] \\
\pp(X,\pp(Z,Y)) \ar[r] & \pp(\pp(X,Z),Y) \ar[r] & \pp(\pp(Z,X),Y)
}$$
commute. 
\hfill $\square$\\
Note that $\theta(p),\theta(q)$ are homotopic for all
$p,q \in \mathcal{P}(2)$. Therefore, the symmetric monoidal structures
$\pp, \underline{q}$ are isomorphic. 
\begin{proposition} \cite[12.2]{GeomILS}\label{RealisationAlgebra}
If $(X_{\bullet},\theta) \in \mathcal{P}[\mathbf{S}\Top]$,
$(|X_{\bullet}|, \tilde{\theta})$ is an algebra over $\mathcal{P}$,
where for $p \in \mathcal{P}$, the action is given by
$\tilde{\theta}(p) := 
\xi \circ |\theta (p)| \circ \xi^{-1}: |X_{\bullet}|^j 
\rightarrow |X_{\bullet}|.$
\end{proposition}
Let $\mathcal{B}(.,.,.)$ 
denote the two-sided bar construction \cite{GeomILS,Markl}.
For a topological group $G$, we use
the models $BG:=|\mathcal{B}(\ast,G,\ast)|$ and 
$EG:=|\mathcal{B}(\ast,G,G)|$ for a classifying space
and a contractible $G$-space. 
\begin{lemma}\cite{GeomILS}\label{POnEG}
If $(G,\theta)$ is a group in $\mathcal{P}[\Top]$,
there are actions of $\mathcal{P}$ on $BG$ and $EG$. Furthermore, the maps
$EG \rightarrow BG$ and $EG \times G \rightarrow EG$ are
 maps of algebras over $\mathcal{P}$.
\end{lemma}
\textbf{Proof:}
First, let us show that $\mathcal{B}(\ast, G, \ast)$ and
$\mathcal{B}(\ast, G, G)$ are algebras over $\mathcal{P}$.
For any $p \in \mathcal{P}(j)$, the map
$\theta_{j} (p): G^j \rightarrow G$ gives rise to maps \\
$\theta_{j} (p)^q : \mathcal{B}_{q}(\ast,G,\ast)^j \cong (G^q)^j \cong (G^j)^q 
\rightarrow G^q \cong \mathcal{B}_{q}(\ast,G,\ast)$
and \\
$\theta_{j} (p)^{q+1} : \mathcal{B}_{q}(\ast,G,G)^j \cong (G^{q+1})^j 
\cong (G^j)^{q+1} 
\rightarrow G^{q+1} \cong \mathcal{B}_{q}(\ast,G,G)$.
Compatibility of $\theta$ with the group structure of $G$ implies that 
these are maps of simplicial spaces. 
Since $$\theta_{j}: \mathcal{P}(j) \rightarrow \Top(G^j,G) \qquad  
\text{and} \qquad 
(.)^{q}: \Top(G^{j},G) \rightarrow \Top((G^{j})^q,G^q)$$
are continuous, 
$\mathcal{B}(\ast,G,\ast)$ and
$\mathcal{B}(\ast,G,G)$ are indeed in $\mathcal{P}[\mathcal{S}\Top]$.
Proposition \ref{RealisationAlgebra} gives $\mathcal{P}$-algebra structures on
$\theta^B$ on $BG$ and $\theta^E$ on $EG$. The 
commutativity of the diagrams
$$
\xymatrix{
\mathcal{B}_q(\ast,G,G)^j \ar[r] \ar[d]^{\theta_{j}^{E}(p)}
 & \mathcal{B}_q(\ast,G,\ast)^j \ar[d]^{\theta_{j}^{B}(p)} \\
\mathcal{B}_q(\ast,G,G) \ar[r] & \mathcal{B}_q(\ast,G,\ast)
}
\hspace{1cm}
\xymatrix{
(\mathcal{B}_q(\ast,G,G)\times G)^j \ar[r] 
\ar[d]^{\theta_{j}^{E}(p) \times \theta_{j}(p)}
 & \mathcal{B}_q(\ast,G,G)^j \ar[d]^{\theta_{j}^{E}(p)} \\
\mathcal{B}_q(\ast,G,G) \times G \ar[r] & \mathcal{B}_q(\ast,G,G)
}
$$ implies that $EG \rightarrow BG$ and $EG \times G \rightarrow EG$
are maps of $\mathcal{P}$-algebras.
\hfill $\square$ \\

If we have an action
$m: G_{+} \wedge F \rightarrow F$ of $G$ on a spectrum
$F \in \spec$, we can define
an ex-spectrum $BGF$ over $BG$ by setting
$BGF(V):=EG \times_{G} F(V)$ with the evident projection and section. 
The $O(V)$ action is the one induced by the action on $F(V)$ and
since $(EG \times_{G} F(V)) \wedge S(W) \cong EG \times_{G} 
(F(V) \wedge S(W))$,
the structure maps of $F$ give rise to structure maps \\
$\sigma_{V,W}: BGF(V) \wedge S(W) \rightarrow BGF(V \oplus W)$.


\begin{lemma}\label{Monoid}
Let $(F,\phi) \in \mathcal{P}[\spec]$, $(G,\theta)$ a group in 
$\mathcal{P}[\Top]$
with an action \\
$m: G_{+} \wedge F \rightarrow F$ in $\mathcal{P}[\spec ]$
and $\mathcal{P}$ a pointed $E_\infty$ operad.
Then for each  $p \in \mathcal{P}(2)$, $BGF$ is a monoid in
$(Ho\spec_{BG},\underline{p})$, i.e. there is a map \\
$\mu(p) : \underline{p}(BGF,  BGF) \rightarrow 
BGF$ and a unit $\eta: 1_{BG} \rightarrow BGF$
satisfying the usual conditions.\end{lemma}
\textbf{Proof:} To construct $\mu(p)$, we have to define
$O(V_1) \times O(V_2)$ -equivariant maps
$$\mu(p)_{V_{1},V_{2}} : 
BGF(V_1) \ewedge BGF(V_2) \rightarrow BGF(V_1 \oplus V_2)$$
for all $V_1,V_2 \in \Ort$.
Recall the construction of the $\mathcal{P}$-actions
$\theta^B$ on $BG$ and $\theta^E$ on $EG$ from Lemma \ref{POnEG}.
The $\mathcal{P}$-structure of $F$ gives $$ 
\phi(p)_{V_1,V_2}: F(V_{1}) \wedge F(V_{2}) \rightarrow F(V_{1}
\oplus V_{2}).$$ Define
\begin{eqnarray*}
\tilde{\mu}(p)_{V_1,V_2}& : & BGF(V_1) \overline{\wedge} 
BGF(V_2)  \rightarrow   \theta^{B}(p)^{\ast} BGF(V_{1}
\oplus V_{2}) \\
 & &  [[e_1,f_1],[e_2,f_2]]  \mapsto 
[\theta^E(p)(e_1,e_2),\phi(p)_{V_1,V_2}(f_1,f_2)].
\end{eqnarray*}
These maps are by Proposition \ref{POnEG} well defined and since they are compatible
with the structure maps as well as $O(V)\times O(W)$ equivariant, we obtain
$$\tilde{\mu}(p): BGF \overline{\wedge}  BGF
\rightarrow \theta^{B}(p)^{\ast} BGF.$$
Letting $\mu(p)'$ denote the adjoint of $\tilde{\mu}(p)$, $\mu(p)$
is given by (the homotopy class of)
$$\pp(BGF,BGF)= \theta^B (p)_!(QBGF \overline{\wedge} QBGF) \rightarrow
\theta^B (p)_!(BGF \overline{\wedge} BGF) \stackrel{\mu(p)'}{\rightarrow} BGF.$$ 
Notice that the same procedure yields maps
$$\mu(q)': \theta^B(q)_!BGF^{\overline{\wedge}^{j}} \rightarrow BGF$$ for each 
$q \in \mathcal{P}(j)$.
Similarly, $\eta: 1_{BG} = i_{b!}S \rightarrow BGF$ is the adjoint of the 
unit $1_F: S \rightarrow F \cong i_{b}^{\ast} BGF$. What remains to be shown
is  the compatibility of $\mu(p)$ and $\eta$ and the 
associativity and unit isomorphisms of $\pp$, i.e. commutativity of the
diagrams
$$
\xymatrix{\pp(BGF,\pp(BGF,BGF)) \ar[rr]^{\sim} \ar[d]^{\pp(id,\mu(p))} & 
&\pp(\pp(BGF,BGF),BGF)  
\ar[d]^{p(\mu(p),id)}  \\
\pp(BGF,BGF) \ar[r]^{\mu(p)}  & BGF
&\pp(BGF,BGF) \ar[l]^{\mu(p)} } $$
and
$$\xymatrix{
\pp(1_{BG},BGF) \ar[r]^{\pp(\eta,id)} \ar[dr]^{\sim}
& \pp(BGF,BGF) \ar[d]^{\mu(p)} \\
& BGF
 }$$

Recall that those isomorphisms were constructed using homotopies given by
paths in the operad. Using the construction above, we define maps of spectra
over these homotopies which we use to prove the result.
 \\
The diagram for the unit transformation
is the following, where $h$ is a path in $\mathcal{P}(1)$ connecting 
$p_{\ast}:=\gamma(p;\ast,id)$ and 1 and $E:=BGF$.\\ 
\pagebreak 
\vspace{1cm}\\

\begin{sideways}
$$\xymatrix{
\\
\theta^{B}(p)_!(Q1_{BG} \ewedge QE) \ar[dr] \ar[d]^{\sim}
\ar[rr]^{\theta^ {B}(p)_!(Q\eta\ewedge id)}
& & \theta^{B}(p)_!(QE \ewedge QE) \ar[dd]& \\
\theta^{B}(p)_!(1_{BG} \ewedge QE) \ar[r] \ar[d]^{\sim}
& \theta^{B}(p)_!(1_{BG} \ewedge E) \ar[dr] \ar[d]^{\sim} 
\ar@{}[ddddrr]^{\mathbf{I}}& & \\
\theta^{B}(p)_! (i_{\ast} \times id)_! (S \ewedge QE) \ar[d]^{\sim} \ar[r] &
\theta^{B}(p)_! (i_{\ast} \times id)_! (S \ewedge E) \ar[d]^{\sim} & 
\theta^{B}(p)_!(E \ewedge E) \ar[dddr]^{\mu(p)} 
& \\
\theta^{B}(p_{\ast})_! QE \ar[r] \ar[d]^{\sim}& 
\theta^{B}(p_{\ast})_! E \ar[ddrr]^{\mu(p_{\ast})}
\ar@{}[ddr]|{\mathbf{II}}
 \ar[d]^{\sim} & & \\
h_! i_{0!}QE \ar[r] \ar[d]^{\sim} &
h_! i_{0!}E \ar[rd] \ar[d] & & \\
h_! \pi^{\ast} RQE \ar[r] & 
h_! \pi^{\ast} RE & 
h_! \pi^{\ast} E \ar[l]^{\sim} \ar[r]^{{\mu}'(h_{t})}& E \\
h_! i_{1!}QE \ar[r] \ar[u]^{\sim} &
h_! i_{1!}E \ar[ru] \ar[u]
\ar@{}[rru]|{\mathbf{III}} & & \\
QE \ar[u]^{\sim} \ar[r]^{\sim} & E \ar[rruu]^{id} \ar[u]^{\sim}
}$$
\end{sideways}
\\
\vspace{2cm}
\\
The commutativity $\mathbf{I,II,III}$ follows from the definiton of
the maps
$\mu(q)'$ for $q \in \mathcal{P}(j)$and the definition of an action of an operad.\\
Similarly, the compatibility with the associativity
transformations is shown by the following diagram.
The top horizontal
line gives the associativity isomorphism and $h$ is the homotopy
corresponding to a path in $\mathcal{P}(3)$ connecting 
$\gamma(p;1,p)$ and $\gamma(p;p,1)$.\\
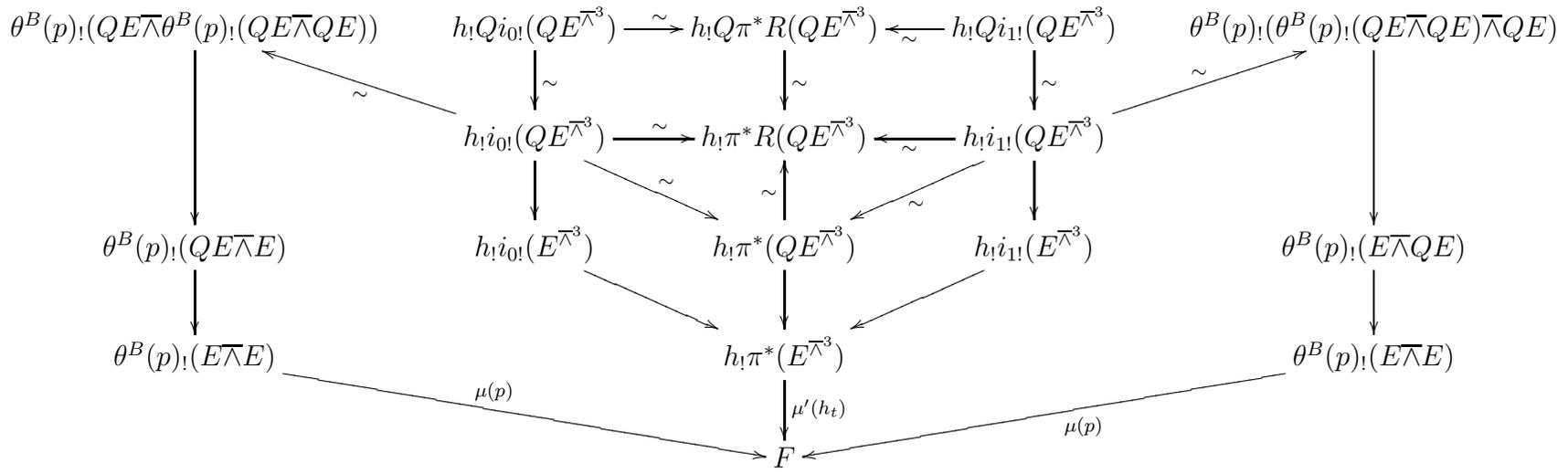
\begin{sideways}
\vspace{3cm}\\
$$\xymatrix{ & & \\
\\
& & & \\
\theta^{B}(p)_{!}(QE \ewedge  	\theta^{B}(p)_{!}(QE \ewedge QE))  \ar[dd] & 
h_! Q i_{0!} (QE^{\ewedge^{3}}) \ar[d]^{\sim} \ar[r]^{\sim} &
h_! Q \pi^{\ast}R(QE^{\ewedge^{3}}) \ar[d]^{\sim} &
h_! Q i_{1!} (QE^{\ewedge^{3}}) \ar[l]^{\sim} \ar[d]^{\sim} &
\theta^{B}(p)_{!}( \theta^{B}(p)_{!}(QE \ewedge QE)\ewedge QE)\ar[dd]
\\
& h_! i_{0!} (QE^{\ewedge^{3}}) \ar[ul]^{\sim} \ar[r]^{\sim} \ar[dr]^{\sim} \ar[d]&
h_!\pi^{\ast}R (QE^{\ewedge^{3}})
&
h_! i_{1!} (QE^{\ewedge^{3}}) \ar[ur]^{\sim} \ar[l]^{\sim} \ar[dl]^{\sim}\ar[d]&
\\
\theta^{B}(p)_{!}(QE \ewedge E) \ar[d] & 
h_!i_{0!}(E^{\ewedge^{3}}) \ar[dr] & 
h_!\pi^{\ast} (QE^{\ewedge^{3}}) \ar[d] \ar[u]^{\sim}&
h_!i_{1!}(E^{\ewedge^{3}}) \ar[dl]
&  \theta^{B}(p)_{!}(E \ewedge QE) \ar[d] \\
 \theta^{B}(p)_{!}(E \ewedge E) \ar[rrd]^{\mu(p)} 
& & h_!\pi^{\ast} (E^{\ewedge^{3}})\ar[d]^{{\mu}'(h_{t})} & & 
 \theta^{B}(p)_{!}(E \ewedge E)\ar[lld]^{\mu(p)}\\
& & F & &  \\
\\
\square
}
$$
\\   
\vspace{2cm} \\
\end{sideways}
\chapter{Twisted cohomology theories}
\section{Definition and basic properties}

Recall that generalised cohomology theories on $\Top_\ast$ 
correspond to spectra, where the theory
corresponding to a spectrum  $E$  can be defined as
$$\tilde{E}^n (X) := [(L\Sigma^{\infty}_{\re^n})X, E] = 
[\Sigma^{\infty}_{\re^n} QX, E] \cong[X, (RE)(\re^n)].$$  By Browns representability theorem, any 
cohomology theory on $\Top_\ast$ arises in this way. 
The situation for twisted cohomology theories is very similar,
where we replace $\Top_\ast$ by the category of ex-spaces $\KB$ and 
$\spec$ by the category of ex-spaces $\specb$.
Just as in the unparametrised case, the morphism sets in $Ho\specb$
can be equipped with an abelian group structure as follows.
$S^2$ is a homotopy commutative H-cogroup with comultiplication given by
the 'pinch map' $S^2 \rightarrow S^2 \vee S^2$. Hence 
$\Sigma^2 E \cong E \wedge S^2$ is an H-cogroup as well for all
$E \in \specb$. This comultiplication defines a commutative product
on $[\Sigma^2 E, F]$ for any $F$ and since $Ho\specb$ is stable,
$$[E,G] \cong [\Sigma E, \Sigma G] \cong [\Sigma^2 E, \Sigma^2 G]$$
is an abelian group for any $E,G \in \specb$. \\ 
We will now define reduced twisted cohomology theories and then use these to define the unreduced 
versions. \\
Recall that in pointed model categories the cofiber of a map
$f: X \rightarrow Y$ is defined to be the coequalizer 
$g: Y \rightarrow Z$ of $f$ and the zero map. If $f$ is a cofibration of 
cofibrant objects, 
$$X \stackrel{f}{\rightarrow} Y \stackrel{g}{\rightarrow} Z$$ is called cofiber sequence.

\begin{definition}
A (reduced) generalised cohomology theory on $\KB$ consists of contravariant 
functors
$\tilde{H}^\ast:  \KB \rightarrow Ab$
indexed on $\mathbb{Z}$ and natural isomorphisms
 $\sigma^n: \tilde{H}^{n + 1} \circ \Sigma \rightarrow
\tilde{H}^{n}$ such that \\
$\bullet$ $ \tilde{H}^\ast$ factors through $\gamma: \KB \rightarrow Ho \KB$; \\
$\bullet$ If $X \rightarrow Y \rightarrow Z$ is a cofiber sequence,
 $\tilde{H}^n (Z) \rightarrow \tilde{H}^n (Y) \rightarrow \tilde{H}^n (X)$ 
is exact for all $n$; \\
$\bullet$ $\tilde{H}^\ast(\coprod_{i \in I} X_i) \cong \prod_{i \in I} 
\tilde{H}^\ast(X_i).$
\end{definition} 

\begin{lemma}
The generalised cohomology theory on $\KB$ correponding to an 
ex-spectrum $E \in \specb$ defined by  
$$\tilde{E}^n (X) := [(L \Sigma_{\re ^{n}}^{\infty}) (X),E] \ \ \ 
\ \ \  f^\ast=\tilde{E}^n(f):= (L\Sigma_{\re^{n}}^{\infty})(f)^\ast \ \  \ 
\forall n \in \mathbb{N}_0
$$
and $$\tilde{E}^{-n} (X) := \tilde{E}^{0}(\Sigma^n X)$$
is indeed a cohomology theory.
\end{lemma}
\textbf{Proof:}
Left derived functors preserve cofiber sequences,
hence 
$$ \tilde{E}^n(Z) \rightarrow \tilde{E}^n(Y) \rightarrow \tilde{E}^n(X) $$ 
is exact for all cofiber sequences $X \rightarrow Y \rightarrow Z$.
For the suspension isomorphism, recall that
the adjoints of the structure maps 
$E(V) \rightarrow \Omega_W E(V \oplus W)$ are weak equivalences
of ex-spaces for all fibrant ex-spectra E, hence
$$[(L \Sigma^{\infty}_{\re^{n+1}}) \Sigma X, E]\cong
[\Sigma X, (R E)(\re^{n+1})] \cong [X, \Omega (R E)(\re^{n+1})] \cong
[X, (R E)(\re^n)].$$

\hfill $\square$

\begin{lemma}\cite{Schoen}
For every generalised cohomology theory $\tilde{H}$,
there is an ex-spectrum $E$ such that
$\tilde{H}^n (X) = [(L\Sigma^{\infty}_{\re ^{n}})X,E] $
\end{lemma}

\textbf{Proof:}
In \cite{Brown}, the author proves a general representability theorem
for functors satisfying the wedge and the Mayer-Vietoris axiom on what he 
calls homotopy categories. 
In \cite[7.5]{MS}, it is shown that $Ho\KB$ is a homotopy category in the sense
of \cite{Brown}, thus the functors $H^n$ 
are representable by ${E}(\re^n)\in \KB$. The transformations
$\sigma^n$ yield the maps $\Sigma {E}(\re^n)\rightarrow E(\re^{n+1}).$
\hfill $\square$\\
To obtain the unreduced theories, we proceed exactly as in the unparametrised case.
\begin{definition}\label{cone}
Let $K_{/B}^{2}$ be the category of pairs in $K_{/B}$. We define the cone
$C: K_{/B}^{2} \rightarrow \KB$ by 
$$C(X,A):=QX \coprod_{QA \times \{0\}} QA\times I \coprod_{QA\times \{1\}} B  \qquad
\text{and}  \qquad C(X,\emptyset):= QX_+.$$ 
\end{definition}
The unreduced cohomology groups are defined by
$$E^n(X,A):=\tilde{E}^n(C(X,A)), \qquad E^n(X):=E^n(X,\emptyset).$$
Note that there are natural maps
$$C(X,A) \rightarrow  C(X,A)/_{B} QX_+ \cong \Sigma QA_+ .$$
Combining these with the suspension isomorphism, we obtain natural maps \\
$$\delta^n: E^{n-1}(A)  \cong \tilde{E}^{n}(\Sigma QA_+) 
\rightarrow \tilde{E}^n (C(X,A)) 
=E^n(X,A).$$ 
\begin{lemma}\cite[5.6.2]{MS}\label{cof}
The sequence
$$..\rightarrow E^{n}(X) \rightarrow E^{n}(A) 
\stackrel{\delta^{n+1}}{\rightarrow} E^{n+1}(X,A) \rightarrow E^{n+1} (X)
\rightarrow ...$$
is exact.
\end{lemma}
The same proof as in the unparametrised case shows that $E^\ast$ satisfies 
excision.
\begin{lemma}\label{exc}
If $X=\stackrel{\circ}{A} \cup \stackrel{\circ}{D}$, the inclusion
$i: (A,A \cap D) \rightarrow (X, D)$
induces an isomorphism
$i^\ast : E^n(X,D) \rightarrow E^n(A, A \cap D).$
\end{lemma}
\textbf{Proof:}
First, let us assume that $X$ is a CW-complex and $A,D$ are subcomplexes.
The natural maps 
$$C(X,D) \rightarrow 
X \cup D\times I \cup B \qquad 
C(A,A\cap D) \rightarrow A \cup (A\cap D) \times I \cup B $$
are weak equivalences by the 5-lemma and \cite[5.6.2]{MS}. 
We use a representation of $(D,A\cap D)$ as NDR-pair to construct a homotopy
inverse of the inclusion
$$\iota: A \cup (A\cap D) \times I \cup B \rightarrow
X \cup D\times I \cup B.$$
Recall that a representation as NDR-pair is given by maps
$$u: D \rightarrow I \qquad h: D \times I \rightarrow D$$
such that
$$u^{-1}(0) = A\cap D , \qquad h_1(u^{-1}[0,1)) = A\cap D, \qquad h_0 = id, 
\qquad 
h_{t|A\cap D} = id.$$
Now, a homotopy inverse of $\iota$ is
$$\begin{array}{lcrclcr}
(A \cup D) \cup & D \times I &\cup B
& \rightarrow & A \cup
& (A\cap D) \times I &
\cup B \\
a & & & \mapsto & a & & \\
\ \ \ \ \ \ \ d & & &\mapsto & & (h_1(d),u(d))  & \\
& (d,t) & & \mapsto & & (h_1(d),max(t,u(d))) & \\
& & b & \mapsto & & & b. 
\end{array}$$
This shows tht $\iota$ is a weak equivalence of ex-spaces and therefore 
induces an isomorphism. To give the proof for general $X$, one uses 
$CW$-substitutes. The crucial observation is that if 
$X=\stackrel{\circ}{A} \cup \stackrel{\circ}{D}$ and one has CW-substitutes
$A'$ of $A$ and $D'$ of $D$ which are obtained as extensions of a CW-substitute 
$(A \cap D)'$ of $A \cap D$, then 
$A'\cup D'$ is a CW-substitute of $X$.
\hfill $\square$  \\
Thus, $E^\ast$ does indeed define a 
cohomology theory on pairs. 
\begin{definition}
A triad $(X; A, D)$ over $B$ is called proper with respect to $E$ if
the inclusions 
$$i_A: (A, A \cap D) \rightarrow (A \cup D, D) \ \  \text{and}\ \ 
i_D: (D, A \cap D) \rightarrow (A \cup D, A)$$ induce isomorphisms
$$i_{A}^{\ast} : E^\ast (A, A \cap D) \rightarrow E^\ast (A \cup D, D) \ \
 \text{and} \ \
i_{D}^{\ast}: E^\ast(D, A \cap D) \rightarrow E^\ast (A \cup D, A).$$
\end{definition}
The Mayer-Vietoris sequence,  
the exact sequence of a triple and the exact sequence of a proper triad follow
formally from Lemma \ref{cof} and \ref{exc}. 
\begin{lemma}
If $A \subset D \subset X$, the sequence
$$...\rightarrow E^{n-1}(A,D) \stackrel{\Delta}{\rightarrow}
E^{n}(X,A) \rightarrow E^n(X,D) \rightarrow E^n(A,D) \rightarrow...$$
where $\Delta$ is the composition $$E^{n-1}(A,D) \rightarrow E^{n-1}(A) 
\stackrel{\delta^{n}}{\rightarrow}E^n(X,A)$$ is exact.
\end{lemma}
\begin{lemma}
If $X=\stackrel{\circ}{A} \cup \stackrel{\circ}{D}$, the sequence
$$...\rightarrow E^{n-1}(A\cup D) \stackrel{\Delta}{\rightarrow}
E^{n}(X) \rightarrow E^n(A)\oplus E^n(D) \rightarrow E^n(A \cup D) \rightarrow...$$
where $\Delta$ is the composition
$$E^{n-1}(A\cup D) \stackrel{\delta^n}{\rightarrow} E^n(A, A\cup D) \cong 
E^n(X,D) \rightarrow E^n(X)$$
is exact. 
\end{lemma}

\begin{lemma}
If $(X;A,D)$ is a proper triad, the sequence
$$...\rightarrow E^{n-1}(A,A\cap D) \stackrel{\Delta}{\rightarrow}
E^{n}(X, A \cup D) \rightarrow E^n(X,D) \rightarrow E^n(A, A \cap D) 
\rightarrow...$$
where $\Delta$ is the composition
$$E^{n-1}(A,A\cap D) \cong E^{n-1}(A \cup D, D) \rightarrow  
E^{n-1}(A \cup D) \stackrel{\delta^n}{\rightarrow} E^n(X,A \cup D)$$
is exact. 
\end{lemma}

There is an additional feature of parametrised cohomology theories
which does not show up in the classical ones, namely the following.
\begin{definition}
If $q,r:X \rightarrow B$ are homotopic via
$h:X \times I \rightarrow B$, the isomorphism induced by $h$ is
$$\psi(h):=i_1^\ast (i_0^\ast)^{-1}:
 E^n(X,A,q) \rightarrow E^n((X,A) \times I,h) \rightarrow  E^n(X,A,r).$$
\end{definition}
In particular, maps $h:X \times S^1 \rightarrow B$ yield
isomorphisms $$\psi(h):E^n((X,A),h_{|X \times \{1\}}) \rightarrow 
E^n((X,A),h_{|X \times \{1\}}).$$  
\begin{lemma}
If two homotopies $h_0,h_1$ are homotopic through
$H$, they induce the same isomorphism.
If $h$ is the composition of the homotopies
$h_0$ and $h_1$, $\psi(h)=\psi(h_1)\psi(h_0)$.
\end{lemma}
\textbf{Proof:}
Consider the diagram
$$\xymatrix{
((X,A),q) \ar[r]^{i_0} & ( (X,A)\times I,h_0)\ar[d]^{j_0} &
((X,A),r) \ar[l]_{i_1} \\
& ((X,A) \times I^2, H) & \\
((X,A),q) \ar[r]^{i_0} & ( (X,A)\times I,h_1)\ar[u]^{j_1} &
((X,A),r). \ar[l]_{i_1} 
 }$$
Observe that  $$j_0 \circ i_0, j_1 \circ i_0: ((X,A),q) \rightarrow 
((X,A)\times I^2, H)$$ as well as 
$$j_0 \circ i_1, j_1 \circ i_1: ((X,A),r) \rightarrow 
((X,A)\times I^2, H)$$
are homotopic over $B$ and therefore 
$$(j_0 \circ i_0)^{\ast} = (j_1 \circ i_0)^\ast \qquad \text{and} \qquad
(j_0 \circ i_1)^\ast = (j_1 \circ i_1)^{\ast}.$$
Thus,
$$\psi(h_0) = i_1^\ast (i_{0}^{\ast})^{-1} =  (j_0 \circ i_1)^{\ast} 
((j_0 \circ i_0)^{\ast})^{-1}=(j_1 \circ i_1)^{\ast} 
((j_1 \circ i_0)^{\ast})^{-1}=\psi(h_1).$$ 
The proof for the second statement is similar.
\hfill $\square$
\section{Product structures}
Recall the monoidal 
structure $\pp$ on $Ho \specb$
for base spaces with an action by an $E_\infty$-operad from Proposition 
\ref{symmstruct}. 
Note that the same construction can be carried through with $\KB$ instead of 
$\specb$. This defines a monoidal structure on $Ho\KB$,
which we again denote by $\pp$.
Observe that there are natural isomorphisms
$$L\Sigma^{\infty}_{V \oplus W} \circ \pp \cong \pp \circ 
(L\Sigma^{\infty}_{V} \times L\Sigma^{\infty}_{W}),$$
justifying the use of the same notation for the monoidal structures on
spaces and spectra.
\begin{definition}
A parametrised ring spectrum over
$B$ is a monoid in $(Ho\specb,\pp)$.
\end{definition}
We use parametrised ring spectra to define products in parametrised 
cohomology theories,
mimicking the use of ring spectra in the construction
of multiplicative theories as in \cite{Adams, Sw}. 
There are three different products,
the wedge product for the reduced theories and the external and internal 
(or cup) product for the unreduced ones.
All of these are very similar to their counterparts in the unparametrised case.
Since the latter are discussed extensively in the literature and the proofs 
carry over to our setting, we will make it short and give little more
than the definitions. 
The formal reason why the proofs carry over is the following.
In \cite{MAddi}, the definition when a triangulation and a closed 
symmetric monoidal structure on a category are compatible is given.
Moreover, he shows how to prove compatibility if the category in question
is the stable homotopy category of a topological model category and the symmetric 
monoidal structure is obtained from a Quillen bifunctor.   
The parametrised stable homotopy category $Ho\specb$ is triangulated
\cite[13.1.5]{MS} and the symmetric monoidal structure $\pp$ on $Ho\specb$
is obtained from a Quillen bifunctor, hence one can follow \cite{MAddi} to proof compatibility.     
\begin{definition}
Let $(E,\mu)$ be a parametrised ring spectrum.
The product 
$$\wedge_p: \tilde{E}^n (X) \times \tilde{E}^m(Y) \rightarrow 
\tilde{E}^{n+m}(\pp(X,Y))$$
in $\tilde{E}$ is given by
\begin{eqnarray*}
\wedge_p: [L\Sigma^{\infty}_{\re^{n}} X, E] \times [L\Sigma^{\infty}_{\re^{m}} Y, E]
& \stackrel{\pp}{\rightarrow} &
[\pp( L\Sigma^{\infty}_{\re^{n}} X,L\Sigma^{\infty}_{\re^{m}} Y),\pp(E,E)]
\\
& \cong & [L\Sigma^{\infty}_{\re^{n+m}} \pp(X,Y),\pp(E,E)]
\\
& \stackrel{\mu_{\ast}}{\rightarrow }&  
[L\Sigma^{\infty}_{\re^{n+m}} \pp(X,Y),E].
\end{eqnarray*}

\end {definition}

The wedge product has the same properties as its non-parametrised analogue,
i.e. it is associative, additive in each variable, 
 commutes with suspension and so on.
\\

Since $\theta(p)_!$ is a left Quillen functor, the natural weak eqivalences
of cofibrant objects
$$C(X,A) \ewedge C(Y,D) \tilde{\rightarrow} C(X \overline{\times} Y, 
X \overline{\times} D \cup A \overline{\times} Y)$$
yield weak equivalences
$$\pp(C(X,A), C(Y,D)) \tilde{\rightarrow} \theta(p)_! C( X \overline{\times} Y, 
X \overline{\times} D \cup A \overline{\times} Y).$$
We abuse notation and write
$(X,A)\overline{\times}(Y,D):=(X\overline{\times}Y, 
A\overline{\times}Y \cup X \overline{\times} D)$.
Using the above equivalences, we obtain the external product
$$\times_p:E^n(X,A) \times E^m (Y,D) \rightarrow E^{n+m} 
\theta(p)_! ((X,A) \overline{\times} (Y,D))$$
and again, this product shares the main features of its non-parametrised
analogue.   
Moreover, the product is compatible with the isomorphisms induced by 
homotopies of twists in the following sense.
If $r,r':Y \rightarrow B$ are homotopic via
$h:Y \times I \rightarrow B$,
$$\theta(p) \circ (q \times h): X \times Y \times I \rightarrow B$$
is a homotopy from $ \theta(p) \circ (q \times r)$  to 
$\theta(p) \circ (q \times r')$
and we have the following commutative diagram:
$$\xymatrix{
E^n((X,A),q) \times E^m((Y,D),r) \ar[r]_{\times_{p}} \ar[d]^{1 \times \psi(h)} &
E^{n+m}(\theta(p)_!((X,A) \overline{\times} (Y,D),q \times r)) 
\ar[d]^{\psi(p \circ (q \times h))} \\
E^n((X,A),q) \times E^m((Y,D),r') \ar[r]_{\times_{p}} &
E^{n+m}(\theta(p)_!((X,A) \overline{\times} (Y,D),q \times r')).}$$ 
Finally, we turn to the parametrised analogue of internal products.
Note that the $\mathcal{P}$-action equips
$[X,B]$ with the structure of an abelian group,
where the sum is given by 
$$[q]+[r]=[q +_p r]:=[p\circ (q \times r) \circ \Delta]$$ 
and the neutral element is the constant map
$pr_0: X \rightarrow \ast \stackrel{i_\ast}{\rightarrow} B$ onto the 
homotopy unit.  
To define internal products, we need parametrised analogues of the diagonal
$$\Delta: (X,A \cup D) \rightarrow (X,A) \times (X,D)$$
and the sum of twists given by the $\mathcal{P}$-action is exactly what is needed
to make this map a morphism over $B$. 
If $q,r:X \rightarrow B$ are two twists, we equip $(X, A \cup D)$ with the twist
$q+_p r$, obtaining  
$$\Delta:((X, A \cup D),q+_p r) \rightarrow
\theta(p)_!((X,A) \times (X,D), q \times r)$$ 
which yields the internal product
$$\cup_p := \Delta^{\ast} \circ \times_p : E^n ((X,A),q) \times 
E^m((X,D),r) \rightarrow E^{n+m} ((X,A \cup D),q +_p r).$$
The $\mathcal{P}$-action can be used to define
associativity and unit homotopies $h_{ass}$ and $h_1$, connecting
$q+_p (r+_p s)$ to $(q+_p r) +_p s$ and
$q+_p pr_o$ to $q$. These are  
unique up to higher homotopies, so induce unique isomorphisms
$$\psi(h_{ass}): E^\ast((X,A),q+_p (r+_p s)) \rightarrow E^\ast((X,A),(q+_p r)+_p s)
$$ and
$$ 
\psi(h_{1}): E^\ast((X,A),q+_p pr_0) \rightarrow E^\ast((X,A),q).$$
We will henceforth make frequent use of these isomorphisms
without further mention. For example, combining the unit isomorphism
and the cup product makes the untwisted groups
$E^\ast((X,A),pr_0)$ a ring. 
Similarly, the twisted groups $E^\ast((X,A),q)$ are $E^\ast((X,A),pr_0)$-modules
for any twist $q: X \rightarrow B$.
\section{Orientation}
Let $\pi: X \rightarrow M$ be a vector bundle. 
We denote the restriction of 
$X$  to $U\subset M$ by $X_U$ and the complement of the zero
section by $X'$. 
\begin{definition}Let $\pi: X \rightarrow M$ be a vector bundle over 
$M \stackrel{q}{\rightarrow} B$ and $\| . \|$ a fibrewise norm on $X$.  
The parametrised Thom space is
$$Th_B (X,q):= X /_B (X \setminus D),$$ 
where $D \subset X$ is the open unit disc
bundle. 
\end{definition}
We will frequently omit the twist from our notation. Since all
bundles are assumed to be bundles over $B$, we use the same
letter for the twist on $M$ and for the corresponding twist on $X$.
Of course, $Th_B(X)$ is weakly equivalent to
the cone $C(X,X')$ which was defined in Definition \ref{cone}.  
From now on, we assume that the base space $B$ is connected,
i.e. that any constant map is homotopic to the constant map
$pr_0$ onto the homotopy unit of $B$.
This implies that
$$E^{i+n}(X_m,X'_m,q_{|m})\cong
\tilde{E}^{i+n}(S^n,pr_{q(m)})\cong \tilde{E}^i(S^0,pr_0)$$
for any $n$-dimensional vector bundle
$X \rightarrow M$ and any $m \in M$. 
In particular,
$E^{\ast}((X_m,X'_m,q_{|m}))$ is a free
$\tilde{E}^\ast(S^0,pr_0)$-module with one generator of degree $n$. 
Note that the above isomorphism depends on the choice of a path 
between $q(m)$ and the homotopy unit. 
\begin{definition}
A vector bundle $X^n \rightarrow M$ over $q:M \rightarrow B$ is called
$(E,q)$-oriented 
by $u \in E^n(X, X') $
if 
$i_{m}^{\ast}u$ is a generator of $E^{\ast}(X_m, X_m')$ as
$\tilde{E}^\ast(S^0,pr_0)$-module.
\end{definition}
\begin{theorem}\label{ThomIso}
If $u$ is an $(E,q)$-orientation of $X^n \rightarrow M$,
and $M$ is a finite CW-complex,
the homomorphism
$$Th(u):E^i(M,r) \stackrel{\pi^\ast}{\rightarrow} E^i(X,r)
 \stackrel{\cup u}{\rightarrow}E^{i+n}((X,X'),r +_p q)
\cong \tilde{E}^{i+n}(Th_B(X,r +_p q))$$ is an isomorphism for  all
 $r:M \rightarrow B$.
\end{theorem}
\textbf{Proof:}
We examine the spectral sequences of the fibrations
$$M \stackrel{id}{\rightarrow} M \qquad \text{and} \qquad
X \rightarrow M$$
and observe that $Th(u)$ induces a morphism of spectral sequences which 
is an isomorphism on the first pages.
Let us first recall the construction of the spectral sequence. 
Denote the preimage of 
the $s$-skeleton 
$M^s$ by
$X^{s}$. 
The exact sequence of the triad
$(X^s;X^{s-1},X'^s)$ is
\begin{eqnarray*}
...& \rightarrow & E^{n-1}(X^{s-1}, X'^{s-1}) \rightarrow
E^{n}(X^{s}, X'^{s}\cup X^{s-1}) 
\\
&\rightarrow&
E^{n}(X^{s}, X'^{s}) \rightarrow
E^{n}(X^{s-1}, X'^{s-1}) \rightarrow
E^{n+1}(X^{s}, X'^{s}\cup X^{s-1}) \rightarrow
...
\end{eqnarray*}
and thus, the groups
$$A^{s,t}:=E^{s+t}(X^s,X'^s), \ \ \ 
C^{s,t}:=E^{s+t}(X^s,X'^s \cup X^{s-1})$$ 
form an exact couple. The standard machinery constructs 
a
spectral sequence with $E_1$-page equal to  $C^{s,t}$
and converging to $E^\ast(X, X')$. \\
For each $s$-cell $j: D^s \rightarrow M$, we denote the bundle
$j^\ast X$ by $X_{D^{s}}$ and regard $D^s$ as a disc over $B$
with the twist given by the pull-back of the twist of $M$.
The map $j$ yields a map of ex-spaces
$$C(X_{D^{s}},X_{S^{s-1}}
\cup X'_{D^{s}})
\rightarrow
C(X^s, X^{s-1}\cup X'^s ).$$
Taking the coproduct over all cells, we obtain
$$
\coprod_{D^s \subset M^s} C(X_{D^{s}},X_{S^{s-1}}
\cup X'_{D^{s}})
\rightarrow
C(X^s, X^{s-1}\cup X'^s ). 
$$
To show that this map is a weak equivalence we use a representation
of $(X^s, X^{s-1})$ as an NDR-pair, i.e. maps
$$u: X^s \rightarrow I \qquad h:X^s \times I \rightarrow X^s$$
such that 
$$u^{-1}(1) = X^{s-1}, \ \ \ h_{t |X^{s-1}} = id, \ \ \  h_0 = id, \ \  \ 
h_1(u^{-1}(0,1]) = 
X^{s-1}.$$ 
We use this representation to define
$$\begin{array}{lcrclcr}
X^s \cup & (X^{s-1} \cup X'^{s}) \times I &\cup B
& \rightarrow & \coprod_{D^s} X_{D^{s}}\cup
& (X_{S^{s-1}} \cup X'_{D^{s}}) \times I &
\cup B \\
x & & &\mapsto & & (h_1(x),u(x))  & \\
& (x,t) & & \mapsto & & (h_1(x),max(t,u(x))) & \\
& & b & \mapsto & & & b. 
\end{array}
$$
Since for $x \in X^s \setminus X^{s-1}$ there is a unique cell
such that $x \in X_{D^s}$ and $X^{s-1} \times I$ is mapped to $B$,
this is welldefined. Note that for general twists, the map is not 
(and cannot be chosen to be) a morphism of ex-spaces. Nonetheless,
it is a homotopy inverse of the map induced by 
$$\coprod_{D^s \subset M^s} (X_{D^{s}}, X_{S^{s-1}} \cup X'_{D^{s}}) 
\rightarrow  (X^s, X^{s-1} \cup X'^{s})$$
which is a morphism of ex-spaces. 
Therefore, we get 
an isomorphism
$$
E^\ast(X^s, X^{s-1}\cup X'^s )
\cong
\bigoplus_{D^s \subset M^s} E^\ast (X_{D^{s}},X_{S^{s-1}}
\cup X'_{D^{s}}).
$$
To analyse the groups
$E^\ast (X_{D^{s}},X_{S^{s-1}}
\cup X'_{D^{s}}),$ 
denote the upper half disc by $D^{s-1}_+ \subset S^{s-1}.$
For any cell
 $D^s \subset M^s$, the exact sequence of the triple
$(X_{D^{s}}, X_{S^{s-1}} \cup X'_{D^{s}}, 
X_{D^{s-1}_{+}}\cup X'_{D^{s}})$
is
\begin{eqnarray*}
...&\rightarrow& E^{n-1}(X_{S^{s-1}}\cup X'_{D^{s}}, 
X_{D^{s-1}_{+}}\cup X'_{D^{s}})) \rightarrow 
E^n(X_{D^{s}}, X_{S^{s-1}}\cup X'_{D^{s}}) 
\\
&\rightarrow&
E^n(X_{D^{s}}, X_{D^{s-1}_{+}}\cup X'_{D^{s}}) \rightarrow 
E^{n}(X_{S^{s-1}}\cup X'_{D^{s}}, 
X_{D^{s-1}_{+}}\cup X'_{D^{s}})\rightarrow...
\end{eqnarray*}
Since the inclusion
$$ X_{D^{s-1}_{+}}\cup X'_{D^{s}} \rightarrow X_{D^{s}}$$
is a weak equivalence,
we obtain isomorphisms 
$$E^{i-1}(X_{S^{s-1}}\cup X'_{D^{s}}, 
X_{D^{s-1}_{+}}\cup X'_{D^{s}}))
\cong E^{i}(X_{D^{s}}, X_{S^{s-1}}\cup X'_{D^{s}}),$$
while excision yields
$$E^{i-1}(X_{S^{s-1}}\cup X'_{D^{s}}, 
X_{D^{s-1}_{+}}\cup X'_{D^{s}}))
\cong
E^{i-1}(X_{D^{s-1}_{-}}, X_{\partial D^{s-1}_{-}}\cup X'_{D^{s-1}_{-}}).
$$
Combining these isomorphisms,
we get
$$E^{i+s}(X_{D^{s}}, X_{S^{s-1}} \cup
X'_{D^{s}}) \cong
E^{i+s-1}(X_{D^{s-1}_{-}}, X_{\partial D^{s-1}_{-}}\cup X'_{D^{s-1}_{-}})\cong...\cong
E^{i}(X_{v_s},X'_{v_{s}})$$ with
$v_s$ being the south pole of $D^s$.
Summarising the above discussion, 
$$C^{s,t}
\cong
\bigoplus_{D^s \subset M^{s}} E^\ast (X_{D^{s}},X_{S^{s-1}}
\cup X'_{D^{s}})\cong
\bigoplus_{D^{s} \subset M^{s}}
E^t(X_{v_{s}},X'_{v_{s}}).$$
Similarly, the spectral sequence of the fibration 
$M \rightarrow M$ is constructed from the exact couple
$$A_{M}^{s,t}:=E^{s+t}(M^s), \qquad 
C_{M}^{s,t}:=E^{s+t}(M^s, M^{s-1}),$$
and we have isomorphisms
$$C^{s,t}_{M} \cong \bigoplus_{D^s \subset M^s} E^t(v_s).$$ 
The restriction of $Th(u)$ to the respective groups
defines a morphism of exact couples.
Now, for any $r: M \rightarrow B$, we have
$$\xymatrix{C^{s,t}_M \ar[rr]^{Th(u)} \ar[d]^{\sim} & &  
C^{s,t}  \ar[d]^{\sim}\\
\bigoplus_{D^s} E^{t}(v_s) \ar[r]^{\pi^{\ast}}&
\bigoplus_{D^s} E^{t}(X_{v_{s}}) \ar[r]^/-.5cm/{\cup i_{v_{s}}^{\ast} u}&
\bigoplus_{D^s} E^{t+n}(X_{v_{s}},X'_{v_{s}}).}
$$  
The bottom horizontal line is an isomorphism by assumption.
The vertical isomorphisms were constructed using only  excision and
the exact sequence of a triple, thus the diagram commutes. 
\hfill$\square$\\
Notice that the only part of the proof where 
we used the finiteness assumption on $M$ was the convergence of 
the spectral sequence. If the
spectral sequence converges for other reasons (for example, if the coefficients of
the cohomology theory are bounded), the theorem is true for arbitrary $M$.\\
The notion of orientation is stable in the sense that if
$u \in E^n(X, X') $ is a $q$-orientation, then
$\sigma^m(u) \in E^{n+m}(X \oplus \re^m,(X \oplus \re^m)')$
is a $q$-orientation as well. 
Sums and pull-backs of oriented bundles are again oriented in the following way. 
If $u \in E^n((X, X'),q)$ is an $(E,q)$-orientation
of a bundle $X \rightarrow M$,
$f^\ast u$ is an $(E,f^\ast q)$-orientation of $f^\ast X$ 
for all $f: N \rightarrow M$.  
If $Y \rightarrow M$ is another bundle over $M$ and
$v \in E^m((Y, Y'),r)$ is its orientation,
$\Delta^\ast (u \times_p v)$ is an $(E,q +_p r)$-orientation of
$X \oplus Y \cong \Delta^\ast (X \times_p Y)$ with the diagonal
$\Delta: M \rightarrow M \times M$. 
\section{Integration}
We will now discuss parametrised push-forward homomorphisms for
oriented maps. 
Let us first recall the construction of push-forwards in the
non-para-metrised setting.
Starting with a smooth proper map
$f:(M,f^\ast q) \rightarrow (N,q)$, we choose a closed embedding
$$\xymatrix{\tilde{f}:M   \ar[rr]^{f \times \iota} 
\ar[drr]^{f}& & N \times \re^m\ar[d]^{\pi_N} \\
& & N}
$$ 
for some large $m$. We will frequently assume that the image
of $\tilde{f}$ is contained in $N \times D^m.$ Note that any two such embeddings
are stably isotopic. This implies that the normal bundle $\nu(\tilde{f})$
is, up to stable isomorphism, independent of the chosen map $\iota$.
The stable isomorphism class of $\nu(\tilde{f})$ is called stable normal bundle.
$f$ is called oriented if the stable normal bundle
of $f$ is oriented; since an orientation class of a bundle defines orientation classes on the stabilisations of the bundle, this notion does make sense.     
Now, let 
$\phi: \nu(\tilde{f}) \stackrel{\sim}{\rightarrow} U \subset N \times \re^m$
 be a tubular neighbourhood. Using $\phi$, one constructs a
collapsing map
$$\hat{f}: N_+ \wedge S^m \rightarrow Th(\nu(\tilde{f}))$$ and defines 
$$f_!: E^\ast(M)\stackrel{Th(u)}{\rightarrow} \tilde{E}^\ast (Th(\nu(f))) 
\stackrel{\hat{f}^\ast}{\rightarrow} \tilde{E}^\ast 
(N_+ \wedge S^m) \cong E^\ast (N).$$ 

We will now mimic the construction of $f_!$ in the parametrised setting.
 
\begin{definition}
A smooth proper map $f:M \rightarrow N$ is called 
$q$-oriented for a twist $q: N \rightarrow B$ 
if (a representative of) the stable normal bundle of $f$ is 
$f^\ast(q)$-oriented.
\end{definition}
Note that we demand that the twist of $\nu(\tilde{f})$ is pulled back from $N$.
The main difficulty lies in the fact that
for general twists $q$, $\phi$ can't be chosen as a 
map over $B$.
However,  the choice of a fibrewise norm $\|.\|$ on 
$\nu(\tilde{f})$
enables us to deform $q:N \times \re^n \rightarrow B$ 
so that at least the restriction of $\phi$ to the unit disc bundle is
indeed a map over $B$. The idea is to make $q$ constant along the 
fibres of the embedded normal bundle.
More precisely,
define $\tilde{q}^{\norm}: N \times \re^m \rightarrow B$ by
$$ \tilde{q}^{\norm} (n,x): =
\left\{ \begin{array}{l c c} 
q(f(m)) &  if & (n,x) = \phi(m,v), \|v\| \leq 1 \\ 
q(\pi_N \phi(m, (\|v\|-1) v) & if &  
(n,x) = \phi(m,v), 1 \leq \|v\|\leq 2 \\
q(n) & else. & 
\end{array}
\right.$$
If we compose $\phi$ with a homeomorphism
from $\nu(\tilde{f})$ to the open unit disc bundle $D\subset\nu(\tilde{f})$, we obtain a tubular neighbourhood over $B$.
Note that we have a homotopy
$h: N \times \re^m \times I \rightarrow N \times \re^m$ such that
$h_1$ is the identitiy and $q \circ h_0 = \tilde{q}^{\norm}$.
We define $\tilde{r}^{\norm}:=r \circ h_0$ for any twist $r: N \rightarrow B$.
Observe that 
$$\widetilde{r+_p r'}^{\norm}= \tilde{r}^{\norm} +_p \tilde{r'}^{\norm}.$$
Write $Y:= N \times \re^m, \ 
Y^{c_{\norm}}:=N \times \re^m \setminus \phi(D)$ and define
the parametrised Pontrjagin-Thom map 
$$\hat{f}^{\norm}: (Y/_B  Y^{c_{\norm}}, \tilde{q}^{\norm}) \rightarrow 
(Th_B(\nu(\tilde{f}), f^\ast q)), \  \ y \mapsto \left\{
\begin{array}{l c}
\tilde{q}^{\norm}(y) & y \notin \phi(D)\\
\phi^{-1}(y) & y \in \phi(D)
\end{array}
\right.
$$
By the construction of $\tilde{q}^{\norm}$, this is a map of ex-spaces.
Now, we define
\begin{eqnarray*}
f_!: E^\ast (M,f^\ast r) & \stackrel{Th(u)}{\rightarrow} & 
\tilde{E}^\ast Th_{B}(\nu(\tilde{f}),f^\ast r +_p f^\ast q)
\stackrel{\hat{f}^{\norm \ast}}{\rightarrow}
\tilde{E}^\ast(Y/_B Y^{c_{\norm}}, \widetilde{r +_p q}^{\norm})
\\
&
\cong & E^\ast((Y, Y^{c_{\norm}}),\widetilde{r +_p q}^{\norm})
\stackrel{\psi(h)}{\rightarrow}
E^\ast((Y,  Y^{c_{\norm}}), r +_p q)
\\
& \rightarrow & 
E^\ast 
(N \times \re^m, N \times (\re^m \setminus D^m)),r +_p q)
\\ &
\cong &
\tilde{E}^\ast(\Sigma^m N_+,r +_p q)
\cong   E^\ast(N,r +_p q).
\end{eqnarray*}
\begin{proposition}
$f_!$ is independent of the choice
of the fibrewise norm, the embedding and the tubular neighbourhood.  
\end{proposition}
\textbf{Proof:}
Let us first show that it is independent of the fibrewise norm.
If $\norm_0,\norm_1$ are two fibrewise norms on $\nu(\tilde{f})$, we can find a path $G$
in the space of norms connecting $\norm_0$ and $\norm_1$.
We use $G$ to show that the morphism  
\begin{eqnarray*}
E^\ast(N \times \re^m,N \times (\re^m \setminus D^m))
& \rightarrow & E^\ast(Y,Y^{c_{\norm_{i}}}) 
\stackrel{\psi(h_{\norm_{i}})}{\rightarrow}
E^\ast(Y,Y^{c_{\norm_{i}}})\\
\rightarrow \tilde{E}^\ast (Y /_B Y^{c_{\norm_{i}}}) 
&\stackrel{(\hat{f}^{\norm_{i}})^\ast}{\rightarrow}& 
\tilde{E}^\ast (Th_{B}^{\norm_{i}}(\nu(\tilde{f}))) 
\stackrel{\sim}{\leftarrow} E^\ast(\nu(\tilde{f}),\nu(\tilde{f})')
\end{eqnarray*}
is the same for $\norm_0$ and $\norm_1$.  
Set $$D^G :=\{(v,t)\ | \ \ \|v\|^{G(t)} < 1\} \subset \nu(\hat{f})\times I$$ and
$$(Y \times I)^{c_{G}}:=N\times \re^m \times I\setminus (\phi \times id_I) (D^G).$$
$\norm_0$ corresponds to the upper path
in the following diagram, whereas the bottom path corresponds to $\norm_1$.
The middle squares commute and the left and right triangle commute up to homotopy, thus the two morphisms agree in $Ho\KB$.
\\
\begin{sideways}
$$\xymatrix{
& & C(Y, Y^{c_{\norm_{0}}}) \ar[r]^{i_0} \ar[d]^{\sim} &
C((Y, Y^{c_{\norm_{0}}})\times I)  \ar[d]^{\sim} & \ar[l]
\\
& C(N \times \re^m, N \times (\re^m \setminus D^m))
\ar[ur] \ar[dr]
& C(Y \times I, (Y\times I)^{c_{G}}) \ar[r] &
 C(Y \times I, (Y\times I)^{c_{G}})\times I) & \ar[l]
\\
& & C(Y, Y^{c_{\norm_{1}}}) \ar[r]^{i_0} \ar[u]^{\sim} &
C(Y, Y^{c_{\norm_{1}}})\times I)  \ar[u]^{\sim}  & \ar[l]\\
& & & & \\
& C(Y, Y^{c_{\norm_{0}}}) \ar[l]^{i_1} \ar[d]^{\sim}
\ar[r]^{\sim} & 
Y/_B Y^{c_{\norm_{0}}} \ar[d]^{\sim}
\ar[r]^{\hat{f}^{\norm_{0}}} &
Th_{B}^{\norm_{0}}(\nu(\hat{f})) \ar[d]^{\sim} & 
\\
 & C(Y \times I, (Y\times I)^{c_{G}}) \ar[l] \ar[r] &
Y \times I /_B (Y\times I)^{c_{G}}) \ar[r] &
\nu(\tilde{f}) \times I /_B (\nu(\tilde{f}) \times I \setminus D^G) &
C(\nu(\tilde{f}),\nu(\tilde{f})') \ar[ul] \ar[dl] \\
&
C(Y, Y^{c_{\norm_{1}}}) \ar[l]^{i_1} \ar[u]^{\sim}
\ar[r]^{\sim} & 
Y/_B Y^{c_{\norm_{1}}} \ar[u]^{\sim}
\ar[r]^{\hat{f}^{\norm_{1}}} &
Th_{B}^{\norm_{1}}(\nu(\hat{f})) \ar[u]^{\sim} &
}
\end{sideways} \\  
The proof that different embeddings and tubular neighbourhoods 
yield the same $f_!$ is similar.
If $$\iota_0:M \rightarrow D^{m_0}, \qquad \iota_1:M \rightarrow D^{m_1}$$ 
are two embeddings and 
$$\phi_0: \nu(\tilde{f}^0) \rightarrow U_0 \subset
N \times D^{m_0}, \qquad \phi_1: \nu(\tilde{f}^1) \rightarrow U_1 \subset
N \times D^{m_1}$$ tubular neighbourhoods, they are (stably) isotopic, i.e.
after stabilising, there are $\iota: M \times I \rightarrow D^m\times I, 
\Phi: \nu(\tilde{f}) \times I \rightarrow U \subset N \times D^{m} \times I$
from $(\iota_0,\phi_0)$ to $(\iota_1,\phi_1)$.
As in the proof for the different norms, these isotopies are used to show that
the different embeddings give the same
$f_!$.
\hfill $\square$\\
Since the push-forward is indpendent of the chosen norm, we will henceforth 
only include the norm in the notation if it is relevant. Otherwise,
we write $\hat{f},Y^c,\tilde{q}, h$ instead of $\hat{f}^{\norm},Y^{c_{\norm}},
\tilde{q}^{\norm},h ^{\norm}$.
\begin{theorem}
If $f$ is $q$-oriented, 
$$f_! (f^\ast (x) \cup y) = x \cup f_! y$$
for all $x \in E^\ast(N,r), \ y \in E^\ast(M,f^{\ast}s)$.
 
\end{theorem}

\textbf{Proof}
The proof is very similar to the one in the unparametrised case, 
see e.g \cite{Dyer}. Let us first show that
$$\hat{f}^{\ast}(\pi_{M}^{\ast}f^\ast x \cup y')= 
\psi(h)^{-1} \pi_{N}^{\ast}x \cup \hat{f}^\ast y' \ \  \forall \ x \in 
\rE(N,r), 
y' \in \tilde{E}^{\ast}(Th_B(\nu(\tilde{f}),f^\ast s)),$$ 
where $h$ is the homotopy needed to make the tubular neighbourhood a map over $B$.
To see this, observe that the following diagram commutes.\\
$$\xymatrix{
 E^\ast(N,r) \times \tilde{E}^{\ast}(Th_B(\nu(\tilde{f}),f^\ast s)) 
\ar[d]^{f^{\ast} \times 1} \ar[r]^{\pi_{N}^{\ast} \times 1}& 
E^\ast(Y,r) \times \tilde{E}(Th_B(\nu(\tilde{f}),f^\ast s))
\ar[d]^{\psi(h)^{-1}\times 1} \\
E^\ast(M,f^\ast r) \times \tilde{E}(Th_B(\nu(\tilde{f}),f^\ast s))
\ar[d]^{\pi_{M}^{\ast} \times 1}
& E^\ast(Y,\tilde{r}) \times \tilde{E}^\ast(Th_B(\nu(\tilde{f}),f^\ast s)) 
\ar[ld]_{\phi^{\ast} \times 1}  \ar[d]^{1 \times \hat{f}^{\ast}} \\
E^\ast(\nu(\tilde{f}),f^\ast r) \times \tilde{E}^\ast(Th_B(\nu(\tilde{f}),f^\ast s))
\ar[d]^{\sim}
& 
E^\ast(Y,\tilde{r}) \times \tilde{E}^\ast((Y,\tilde{s})/_B Y^{c}) 
\ar[d]^\sim 
\\
E^\ast(\nu(\tilde{f}),f^\ast r) \times E^\ast((\nu(\tilde{f}),
\nu(\tilde{f})\setminus D),f^\ast s) \ar[d]^{\cup}
&
E^\ast(Y,\tilde{r}) \times E^\ast((Y,Y^{c}),\tilde{s} )
\ar[l]_{ \ \  \ \  \  \ \  \ \ \ \phi^{\ast} \times \phi^{\ast}} 
\ar[d]^{\cup}
\\
E^\ast((\nu(\tilde{f}),
\nu(\tilde{f})\setminus D),f^\ast (r +_p s)  \ar[d]^\sim 
&
E^\ast((Y,Y^{c_{g}}),\tilde{r} +_p \tilde{s} )\ar[d]^\sim
\\
E^\ast(Th_B(\nu(\tilde{f}),f^\ast(r +_p s)))
\ar[r]^{\hat{f}^{\ast}}
& \tilde{E}^\ast((Y,\tilde{r} +_p \tilde{s})/_B Y^{c_{g}}) 
}$$
Starting with $x \in \rE(N,r), 
y' \in \tilde{E}^{\ast}(Th_B(\nu(\tilde{f}),f^\ast s))$, following the arrows 
on the left it is mapped to $\hat{f}^{\ast}(\pi_{M}^{\ast}f^\ast x \cup y')$
whereas it is mapped to   \\
$\psi(h)^{-1} \pi_{N}^{\ast}x \cup \hat{f}^\ast y'$ by the arrows
on the right,
proving the equality
$$\hat{f}^{\ast}(\pi_{M}^{\ast}f^\ast x \cup y')= 
\psi(h)^{-1} \pi_{N}^{\ast}x \cup \hat{f}^\ast y'.$$
Recall that 
$$f_!(x'):= \sigma^{-n} i^\ast 
\psi(h)\hat{f}^\ast(\pi_{M}^{\ast}x \cup u)$$
with $i:(N\times\re^m, N\times(\re^m-D^m)) \rightarrow
(Y,Y^{c_{g}})$. Using the above equality,
\begin{eqnarray*}
f_!(f^{\ast}x \cup y) & = &
\sigma^{-n}i^\ast \psi(h) 
\hat{f}^{\ast}(\pi_{M}^{\ast}(f^{\ast}x \cup y)\cup u) \\
&=& 
\sigma^{-n}i^\ast \psi(h) \hat{f}^{\ast}( 
\pi_{M}^{\ast} f^{\ast}x \cup 
(\pi_{M}^{\ast}y\cup u)) \\
&=&
\sigma^{-n}(\pi_{N}^{\ast}x \cup 
i^\ast \psi(h)\hat{f}^{\ast}(\pi_{M}^{\ast} y \cup u)) \\
&=& x \cup f_! y.
\end{eqnarray*}
\hfill $\square$\\

Let us now discuss the functoriality of the push-forward,
so take 
$$f:M \rightarrow N,\ g: N \rightarrow L$$ smooth and assume
that $f$ is $g^\ast q$-oriented  by 
$u \in E^\ast((\nu(\tilde{f}),\nu(\tilde{f})'),f^\ast g^\ast q)$
and $g$ is $r$-oriented by 
$v\in E^\ast((\nu(\tilde{g}),\nu(\tilde{g})'), g^\ast r)$.
We choose 
tubular neighbourhoods
$$\phi^M: \nu(\tilde{f}) \rightarrow N \times \re^m \qquad \qquad 
\phi^N: \nu(\tilde{g}) \rightarrow L \times \re^n$$
and fibrewise norms on $\nu(\tilde{f}), \nu(\tilde{g}).$
To see that there is an induced orientation of $g \circ f$,
observe that
the normal bundle of
$$\widetilde{gf}:=(\tilde{g}\times 1) \circ \tilde{f}
:M \rightarrow L \times \re^{n+m}$$  is
$$\nu(\widetilde{gf}) \cong 
\nu(\tilde{f}) \oplus \tilde{f}^\ast \nu(\tilde{g}\times 1) \cong
\nu(\tilde{f}) \oplus f^\ast \nu(\tilde{g}).
$$
We denote the induced orientation by 
$$uv:=\Delta^\ast (u \times_p f^\ast v) \in E^\ast((\nu(\widetilde{gf}),
\nu(\widetilde{gf})'),f^\ast g^\ast (q +_p r)).$$
Moreover, the tubular neighbourhoods of
$\tilde{f}$ and $\tilde{g}$ can be used to construct a tubular 
neighbourhood of $\widetilde{gf}$ as follows.
For any $w \in \nu(\tilde{f})_{m}$, there is a natural
path $\tilde{\gamma}_w$ joining $\tilde{f}(m)$ and $\phi^M(w)$, namely
$$\tilde{\gamma}_w: I \rightarrow N \times \re^m, \ t \mapsto \phi^M(tw).$$
Choose a connection and thus the notion of parallel
transport on $\nu(\tilde{g}\times 1) \cong \nu(\tilde{g})\times \re^n$
which is compatible with the chosen norm on $\nu(\tilde{g})$ . 
We denote the isomorphism given by parallel transport
along $\tilde{\gamma}_w$ by
 $$\gamma_{w}:  \nu(\tilde{g}\times 1)_{|\tilde{f}(m)} \rightarrow 
\nu(\tilde{g}\times 1)_{|\phi^m(w)}.$$  
A tubular neighbourhood is
$$\phi:=(\phi^N\times 1)\circ \phi^1:
\nu(\widetilde{gf}) \rightarrow L \times \re^{n+m}$$
with\\
 $$
\begin{array}{r c r}
\phi^1:\nu(\widetilde{gf}) \cong 
 \nu(\tilde{f}) + \tilde{f}^\ast \nu(\tilde{g}\times 1) & \rightarrow &
\nu(\tilde{g})\times \re^m \cong  \nu(\tilde{g}\times 1)\\
(w_1, w_2) & \mapsto&  \gamma_{w_{1}} (w_2).
\end{array} 
$$ 
\begin{theorem}
$(g\circ f)_! = g_! \circ f_!$ for the induced orientation of $gf$.
\end{theorem}
\textbf{Proof:} 
We use the tubular neighbourhood $\phi$ defined above and the norm on 
$$\nu(\tilde{gf})=\nu(\tilde{f})+f^\ast \nu(\tilde{g})  \qquad 
\text{given by} \qquad 
\|w_1 + w_2\|:=max\{\|w_1\|,\|w_2\|\}.$$
The advantage of this norm is that 
the various adjustments of the twists needed to
make the tubular neighbourhoods maps over $B$  and the homotopies between
the original twists and the deformed ones are compatible in the following sense. 
The homotopy
$$h^M: N \times \re^m \times I \rightarrow N \times \re^m$$ making
$\phi^M$ a map over $B$ is just the restriction of the
homotopy $h^1$
which turns $\phi^1$ into a map over $B$. 
We abuse notation and 
denote all homotopies by $h$.
Since $(g\circ f)_!$ depends only on the orientation, it is sufficient to
prove the equality for these choices. 
We write $D$ for the unit disc bundle of a Riemannian vector bundle and
$$\begin{array}{ccc}
Y:= N \times \re^m  & Y^c:= N \times \re^m \setminus \phi^M(D) &
Z:= \nu(\tilde{g}\times 1) \\ 
Z^c:=\nu(\tilde{g}\times 1) \setminus \phi^1(D) &
W:=L \times \re^{n+m} & W^c:= L \times \re^{n+m} \setminus \phi(D)\\
X:=L \times \re^n & X^c:=L \times \re^n \setminus \phi^N(D)\\
TF:=Th_B(\nu(\tilde{f})) & TG :=Th_B(\nu(\tilde{g})) &
TGF:=Th_B(\nu(\widetilde{gf})).
\end{array}$$ 
$\hat{gf}:W /_B W^c
\stackrel{G}{\rightarrow}Z /_B Z^c \stackrel{F}{\rightarrow} 
TGF$ factors over
$Z /_B Z^c$, using the factorisation of $\phi=(\phi^N \times 1) \circ \phi^1$.  
In the following commutative diagramm, the left hand vertical morphism is 
$f_!$, the bottom horizontal one is $g_!$ and the diagonal is
$(g\circ f)!.$ \\
\begin{sideways}
$$\xymatrix{
E^\ast(M) \ar[ddrr]^{Th(uv)} \ar[dd]^{Th(u)} 
& & & & &  \\
&
& & & &  \\
\rE(TF)
\ar[r]^/-.5cm/{\wedge f^{\ast} v} \ar[d]^{\hat{f}^{\ast}}
& 
\rE \pp(TF,
Th_B f^{\ast}\nu(\tilde{g})) 
\ar[r]^{\Delta^{\ast}}
&
\rE(TGF)
\ar[d]^{F^\ast} \ar[dr]^{\hat{gf}^{\ast}}
& & & 
\\
\rE(Y /_B Y^{c} 
) 
\ar[r]^{\wedge_{p} v} \ar[d]^{\psi(h)}
&
\rE(\pp(Y/_B Y^{c},
TG)) \ar[r]^{\Delta^{\ast}}
&
\rE(Z /_B Z^c)
 \ar[d]^{\psi(h)} \ar[r]^{G^{\ast}}
& \rE(W/_B W^c) \ar[dr]^{\psi(h)}
& & 
\\
\rE(Y/_B Y^{c}) 
 \ar[d]
& & 
\rE(Z /_B Z^c
 ) \ar[d]
& & \rE(W /_B W^c) \ar[d]  &   
\\
\rE(Th_B(Y)) \ar[r]^{\wedge_p v} \ar[d]^{\sim}
&\rE\pp(Th_B(Y),
TG) \ar[r]
  &
\rE(Th_B(\nu(\tilde{g}) + \re^m)) \ar[d]^{\sim}   
& &  \rE(\Sigma^{n+m}L_+) \ar[dl]^{\sigma^{-m}} \ar[dd]^{\sigma^{-m-n}} 
\\
\rE(\Sigma^m N_+) 
\ar[d]^{\sigma^{-m}}
& &
\rE(\Sigma^m TG) \ar[d]^{\sigma^{-m}}
&  
\rE(\Sigma^n L_+) \ar[dr]^{\sigma^{-m}} &
\\
E^\ast(N) \ar[rr]^{Th(v)}
&
&
\rE( TG)\ar[r]^{\psi(h) \circ\hat{g}^{\ast}} 
& \rE(X /_B X^c) \ar[u] &
E^\ast(L)
\\
\square
}$$
\end{sideways}
\chapter{Examples}
\section{Twisted K-theory}
We will now define $BPU$-twisted (complex) $K$-theory.
In \cite{BJS}, the authors developed
an orthogonal (strict) ring spectrum representing 
$K$-theory. Their construction uses the picture
of $KK$-theory from \cite{BJ}, where unbounded operators are used
to represent $KK$-classes. 
We modify the construction of the $K$-spectrum
in the sense that we weaken the product (obtaining an $E_\infty$-spectrum),
but make way 
for a $PU$-action on the spectrum
which is compatible with the product structure.
First, we recall the results from
\cite{BJS}.
\\
Let $A$ be a $C^{\ast}$-algebra and $M$ a countably generated 
$\mathbb{Z}_2$-graded Hilbert
$A$-module. 
We denote by $B(M)$ the $C^{\ast}$-algebra of bounded $A$-linear  
operators that admit an adjoint with respect to the $A$-valued scalar
product.

\begin{definition}
$ \FM :=\{F \in B(M)\ | \ F^{\ast} = F,\  F^2-1\  compact\ and \ ||F|| \leq 1 \}$ \\
We equip $\FM$ with the weakest topology such that
the maps
$$\FM \rightarrow B(M), \ F \mapsto F^2 \qquad
\text{and} \qquad  
\FM \rightarrow M, \ F \mapsto F(\psi)$$
are continuous (with respect to the norm topology on $B(M)$ and $M$) 
for all $\psi \in M$.
\end{definition}

Let $\textbf{U}_M$ be the subspace  of unitary elements in $\FM$
and denote by $\FM / \textbf{U}_M$ the quotient (in topological spaces).

\begin{definition} 
 $\FM^+$ is the set $\FM / \textbf{U}_M$ with the topology
generated by the open sets in $\FM \setminus \textbf{U}_M$,
and for all $\varepsilon > 0$ the sets
$$\{F \in \FM \, | \, 0 < || F^2-1|| < \varepsilon\} \cup \ast.$$
\end{definition} 

\begin{proposition}
The conjugation action of the unitary group
$U(M)\subset B(M)$ (with the norm topology) 
on $\FM$ is continuous and extends to a
continuous action on $\FM^+$.
\end{proposition}

\textbf{Proof:}
Let us first show that $U(M)$ acts continuously on $\FM$.
Observe that for $U,U_0 \in U(M)$  we have
$\| U^{-1} - U^{-1}_{0} \| = \| U_0 (U_{0} - U) U^{-1}_{0} \| =
\| U - U_{0} \|.$  
\\ 
Now, for $\psi \in M$ and $ F,F_0 \in \FM$ we have
\begin{eqnarray*}
\| UF U^{-1} \psi - U_{0}F_{0} U_{0}^{-1} \psi \| & \leq &
 \| UF U^{-1} \psi - U F U_{0}^{-1} \psi \|  \\ & & + 
\| UF U^{-1}_{0} \psi - U F_{0} U_{0}^{-1} \psi \| \\ & & +
\| UF_0 U^{-1}_{0} \psi - U_{0}F_{0} U_{0}^{-1} \psi \| \\
& \leq &
2 \|U - U_0 \| \| \psi \| + \|(F-F_0) U_{0} \psi \|.
\end{eqnarray*}
Moreover,
\begin{eqnarray*}
\|(U F U^{-1})^2 - (U_0 F_0 U^{-1}_{0})^2 \| & \leq & 
\|  U F^2 U^{-1} -U F^2 U_{0}^{-1} \| \\
& & + \| U F^2 U_{0}^{-1} - 
U F^{2}_{0} U_{0}^{-1} \| 
\\ & & + 
\| U F^{2}_{0} U_{0}^{-1} - U_{0} F^{2}_{0} U_{0}^{-1} \| \\
& \leq & 2 \|U - U_{0} \| + \|F^2 -F_{0}^{2} \| \\
\end{eqnarray*} 
so the action of $U(M)$ on $\textbf{F}_M$ is indeed continuous. 
Since $\textbf{U}_M$ is invariant under this action and
$\|F^2-1 \| = \|U F^2 U^{-1} -1\|$, it extends to an action on $\FM^+$. 
\hfill $\square$ \\

Recall that an $A$-linear operator $D$ on $M$ is called regular if it is
closed densely defined, dom $D^\ast$ is dense and $(1+D^\ast D)$
has dense image. The unbounded picture of $KK$-theory as developed in \cite{BJ}
has the advantage that it allows for an easier description of the 
product. 

\begin{definition} \cite{BJS} \\
$$\RM:=\{D \ regular \ on \ M \, | \,  D^\ast=D, 
\ D \ odd, (1+D^2)^{-1} \ compact\}. $$
\end{definition}
Let 
$J: \RM \rightarrow \FM$ be the map $D \mapsto D (1+D^2)^{-\frac{1}{2}}$.
$J$ 
maps $\RM$  indeed to $\FM$,  is injective, 
and $im(J) \cap \textbf{U}_M = \emptyset$,
hence we may regard $\RM$ as a subset of $\FM^+$\cite{BJ,BJS}. 
\begin{definition}
Let $\RM^+$ be the pointed set $\RM \cup \ast$, equipped with the subspace
topology induced by $\RM^+ \subset \FM^+$.
\end{definition}
Since $\RM^+ \subset \FM^+$ is invariant under the
conjugation action of the even part $U^{0}(M)$ of $U(M)$, the $U(M)$-action
on $\FM$ induces a continuous 
action of $U^{0}(M)$ on $\RM^+$. \\
If $M$ and  $M'$ are Hilbert modules over $A$ and $A'$, we define 
$$\mu: \RM \times \mathbf{R}_{M'} \rightarrow \mathbf{R}_{M \otimes M'},
\ \ (D,D') \mapsto D \otimes 1 + 1 \otimes D'.$$
\begin{proposition}\cite[3.6]{BJS}
$\mu: \RM \times \mathbf{R}_{M'} \rightarrow \mathbf{R}_{M \otimes M'}$ is
continuous and has a continuous pointed extension
$\mu^+: \RM^+ \wedge \mathbf{R}_{M'}^{+} \rightarrow 
\mathbf{R}_{M \otimes M'}^{+}$. 
\end{proposition}

For a finite dimensional real inner product space $(V,g)$, we denote the
Clifford algebra of 
$(V \otimes \mathbb{C}, g \otimes \mathbb{C})$ by $\cx l(V)$.
Note that $\cx l(V)$ is functorial in $(V,g)$ and that 
there are natural isomorphisms
$\cx l(V) \otimes \cx l(W) \cong \cx l(V\oplus W)$.
Consider the $\mathbb{Z}_2$-graded vector space
$C_{0}^{\infty}(V,\cx l(V))$ and define $E_V:=L^2(V,\cx l(V))$ to be its $L^2$-completion.
\begin{definition}\cite{BJS}
\begin{itemize}
\item The Dirac operator $\partial_V \in \mathbf{R}_{E_V}$ is the closure of
the operator on \\ 
$C_{0}^{\infty}(V,\cx l(V))$ given on the functions of pure degree by
$$\partial_V \sigma (x):= (-1)^{|\sigma|}\sum_{i \in I} \left(\frac{\partial \sigma}{\partial v_i}(x)\right) v_i, \qquad x \in V$$
where $(v_i)_{i \in I}$ is any orthonormal basis of $V$.
\item The Clifford operator $L_V \in  \mathbf{R}_{E_V}$ is defined by
$$L_V (\sigma) (v) = v \sigma(v), \ \ \ v \in V, \ \ \ \sigma \in C_{0}^{\infty}(V,\cx l(V))$$ 
\item The Bott-Dirac operator $D_V \in \mathbf{R}_{E_V}$ is 
$D_V := \partial_V + L_V$.
\end{itemize}
\end{definition}
\begin{proposition}\cite[3.9]{BJS}
There are canonical isomorphisms \\
$E_V \otimes E_W \cong E_{V \oplus W}$.
Using these isomorphisms, we have identities
$$\partial_{V \oplus W} = \mu(\partial_V, \partial_W) \qquad \text{and} \qquad 
D_{V \oplus W} = \mu (D_V, D_W).$$ 
 \end{proposition}

Recall that 
an orthogonal spectrum is defined as a module
over the orthogonal sphere spectrum, 
therefore the data for an orthogonal spectrum $K$
is given by 
\begin{itemize}
\item a pointed space $K(V)$ with a left $O(V)$-action
\item $O(V) \times O(W)$-equivariant maps 
$\sigma_{V,W}: K(V) \wedge S^W \rightarrow K(V \oplus W)$
\end{itemize}
for all finite dimensional real inner product spaces $V,W$.
Moreover, 
the appropriate associativity and unit diagrams have to commute. \\
To obtain the spaces of the spectrum representing complex $K$-theory,
we regard $H'(V):=\cx l(V) \otimes E_V$  as a right $\cx l(V)$-module,
where $\cx l(V)$ acts by right multiplication on $\cx l(V)$ and trivially on $E_V$,
and define
$K'(V):=\mathbf{R}_{H'(V)}^{+}$.
\\
To construct the $O(V)$-action, let
$Pin^{c}(V) \subset \cx l(V)^{\ast}$ be the subgroup generated
by elements of $V$ of norm one. 
The $Pin^c$-representation $$\rho_V: Pin^{c}(V) \rightarrow O(V)$$
is defined as follows.
The subspace $V \subset \cx l(V)$ is invariant under the twisted conjugation action 
$$Pin^{c}(V) \times \cx l(V) \rightarrow \cx l(V), \qquad
(g,v) \mapsto (-1)^{|g|}g v g^{-1} $$
and $\rho_V$ is simply the restriction of this action to $V$.
$Pin^{c}(V)$ is a central $S^1$-extension of $O(V)$, i.e.
$$1 \rightarrow S^1 \rightarrow Pin^{c}(V) \stackrel{\rho_V}{\rightarrow} O(V) \rightarrow 1.$$

We define an action of $Pin^{c}(V)$ on $E_V$ by
$$g\sigma(v):= \rho_V(g) \sigma(\rho_V(g)^{-1} v)$$ for
$g \in Pin^{c}(V), \sigma \in C_{0}^{\infty}(V,\cx l(V))$ and $v \in V$. 
Combining this with the $Pin^{c}(V)$-action on $\cx l(V)$ given by left multiplication,
we obtain a unitary representation
$\varrho'_V: Pin^{c}(V) \rightarrow \mathbf{U}_{H'(V)}$.
The twisted conjugation action on $\mathbf{F}_{H'(V)}$ given by
$p: F \mapsto (-1)^{deg(p)} \varrho'_V(p) F \varrho'_V(p)^{-1}$ 
factors through an $O(V)$-action. Moreover,
this $O(V)$-action induces one on $\mathbf{F}_{H'(V)}^{+}$ \cite{BJS}.
The subspace $\mathbf{R}_{H'(V)}^{+}$ is $O(V)$-invariant, 
thus we get an action on 
$K'(V)$.\\   
For $v \in V$, denote the operator on $\cx l(V)$ given by left Clifford multiplication by 
$l_v$ and observe that the map
$$\tilde{\eta}_V: V \rightarrow \mathbf{R}_{H'(V)}, \ \ v \mapsto \mu(l_v, D_V)$$
has a continuous pointed extension 
$\eta_V: S^V \rightarrow \mathbf{R}_{H'(V)}^+=K'(V)$
since 
$$\|J(\tilde{\eta}_{V}(v))^2 -1 \|=\| \frac{(l_v \otimes 1 + 1 \otimes D_V)^2}
{1 + (l_v \otimes 1 + 1 \otimes D_V)^2} - 1 \| 
\stackrel{v \rightarrow \infty}
{\longrightarrow} 0.$$
We define the structure maps of $K'$ by
$$K'(V) \wedge S^W \stackrel{1\wedge \eta_W} {\longrightarrow} 
K'(V)\wedge K'(W) \stackrel{\mu^+} {\longrightarrow}
\mathbf{R}^{+}_{H'(V)\otimes H'(W)}
\cong \mathbf{R}^{+}_{H'(V\oplus W)}=K'(V\oplus W).$$
These maps are $O(V) \times O(W)$-equivariant, associative and unital,
hence $K'$ is indeed an orthogonal spectrum.

\begin{proposition}\cite{BJS}
$K'$ represents complex K-theory. \\
Moreover,
$$\mu: K' \wedge K' \rightarrow K' \ \text{given by} \ 
K'(V) \wedge K'(W) \stackrel{\mu^+}{\rightarrow}
 \mathbf{R}^{+}_{H'(V\oplus W)}=K'(V\oplus W)$$
induces the product in $K$-theory and $\mu,\eta$ make $K'$ a ring spectrum.
\end{proposition}
With this definition of the $K$-spectrum at hand, we will now describe the 
modification that enables us to define $BPU$-twisted K-theory. 
The spaces of the spectrum $K'$ are given as spaces of operators on
Hilbert modules indexed on $\Ort$. The idea is to tensor all these
modules by a fixed Hilbert space insomuch that the projective
unitary group of this space acts on the spectrum, so let us
fix the separable  $\mathbb{Z}_2$-graded Hilbert space 
$\mathbf{H}:=L^2(\re,\cx l(\re))$ 
and define $$
PU:= U^{0}(\mathbf{H}) / S^1 \ \ \ \ \ \ 
H(V):=H'(V) \otimes \mathbf{H} \ \ \ \ \ \ K(V):=\mathbf{R}^{+}_{H(V)}.$$
To get the $O(V)$-action, we extend the unitary $\cx l(V)$-representation \\
$\varrho'(V): Pin^{c}(V) \rightarrow \mathbf{U}_{H'(V)}$ trivially, i.e.
$$\varrho(V):=\varrho'(V) \otimes 1 : Pin^{c}(V) \rightarrow \mathbf{U}_{H(V)}=
\mathbf{U}_{H'(V)\otimes \mathbf{H}}.$$
As before, the twisted conjugation action
factors through an $O(V)$-action, which yields an $O(V)$-action on $K(V)$.
Similarly, we define the structure maps by
$$K(V) \wedge S^W \stackrel{1\wedge \eta_W} {\longrightarrow} 
K(V)\wedge K'(W) \stackrel{\mu^+} {\longrightarrow}\mathbf{R}_{H(V)\otimes H'(W)}
\cong \mathbf{R}_{H(V\oplus W)}=K(V\oplus W).$$
The $O(V) \times O(W)$-equivariance as well as the unitality and associativity 
follow as before, so we have an orthogonal spectrum, again.
The $PU$-action defined by
$$PU \times K(V) \ \ \ \ \ \  
([u],D) \mapsto (1 \otimes u) D (1 \otimes u)^{-1}$$ is
$O(V)$-equivariant as well as compatible with the structure maps and
thus yields a $PU$-action on $K$. 

Moreover, 
$$\phi_V := \mu(., D_{\re}): K'(V) = \mathbf{R}^{+}_{H'(V)} \rightarrow 
\mathbf{R}^{+}_{H(V)}=K(V)$$ 
is a weak equivalence for all $V \neq \{0\}$ \cite[3.12]{BJS}. 
Since these maps are $O(V)$-equivariant and compatible with the 
structure maps, they assemble to a stable equivalence
$\phi: K' \rightarrow K$.

\begin{definition}
Let $\mathcal{P}$ be the $E_{\infty}$-operad given by
$\mathcal{P}(j):= \textup{Isom}^{0}
(\mathbf{H}^{ \otimes j}, \mathbf{H})$
with the obvious $\Sigma_j$-action and structure maps,
where $\textup{Isom}^{0}$ denotes the space of even isometric isomorphisms
with the norm topology. 
\end {definition}
Recall that $\mathcal{P}[\Top]$ is the set of algebras over $\mathcal{P}$,
i.e. spaces $X \in \Top$ with morphisms
$\theta_j: \mathcal{P}(j) \times X^j \rightarrow X$ which are compatible
with the operadic structure of $\mathcal{P}$.
\begin{lemma}
$PU$ is a group object in  $\mathcal{P}[\Top]$.
\end{lemma}
\textbf{Proof:}
Define the action of $\mathcal{P}$ by
$$\theta: \Pj \times PU^j \rightarrow PU, \ \ \ 
(p,[u_1],..,[u_j]) \mapsto [ p  (u_1 \otimes . . . \otimes u_j ) p^{-1}]$$
and check that it commutes with the product
$\mu:PU \times PU \rightarrow PU$ and the unit
$e:\ast \rightarrow PU$. 
\hfill $\square$\\
Since the category of spectra $\spec$ is enriched, tensored and cotensored over
$\Top$, the category $\mathcal{P}[\spec]$ of $\mathcal{P}$-algebras in $\spec$
is defined.
\begin{lemma}
$K$ is in $\mathcal{P}[\spec]$.
\end{lemma}
\textbf{Proof:}
Each $p\in \Pj$ defines a map $\textbf{R}_{M \otimes 
\textbf{H} ... \otimes \textbf{H}} \rightarrow \textbf{R}_{M \otimes 
\textbf{H}}$ for all $M$, given by
$D \mapsto (1 \otimes p) D (1 \otimes p)^{-1}$. 
Using this, the action of $\mathcal{P}$ is defined by the maps 
\begin{eqnarray*}
\phi_{V_{1},..,V_{j}}: \Pj_+ \wedge K(V_1) \wedge ... \wedge K(V_j) 
& \stackrel{1 \times \mu^+}{\rightarrow}& 
\Pj \times \textbf{R}_{H'(V_1 \oplus  ... \oplus V_j) \otimes 
\textbf{H} ... \otimes \textbf{H}} \\
& \rightarrow &  \textbf{R}_{H'(V_1 \oplus  ... \oplus V_j) \otimes 
\textbf{H}}\\
&= &K(V_1 \oplus...\oplus V_j).
\end{eqnarray*}
\hfill $\square$ \\

\begin{lemma}
The $PU$-action on $K$ is compatible with the 
$\mathcal{P}$-actions on $PU$ and $K$, i.e.
$m: PU_+ \wedge K \rightarrow K$ is a map in  
$\mathcal{P}[\spec]$.
\end{lemma}
\textbf{Proof:}
We have to show commutativity of 
$$\xymatrix{
\mathcal{P}(j)_+ \wedge PU_+ \wedge K(V_1)\wedge ... \wedge PU_+ \wedge K(V_j) 
\ar[r] \ar[d]  & 
\mathcal{P}(j)_+ \wedge  K(V_1)\wedge ... \wedge K(V_j) \ar[d] \\
PU_+ \wedge K(V_1 \oplus ... \oplus V_j) \ar[r] & K(V_1 \oplus ... \oplus V_j).
}$$
On inspection of the definitions, this reduces to the observation
that for all $D_i \in K(V_i)$ and $[u_i] \in PU$ \\
$\mu^+((1\otimes u_1)D_1 (1\otimes u_1)^{-1},...,
(1\otimes u_j)D_j (1\otimes u_j)^{-1}) =  \\ 
 \hspace*{\fill}  (1\otimes u_1 ...\otimes u_j)\mu^+(D_1,..,D_j) 
(1 \otimes u_1 ...\otimes u_j)^{-1}.\square$ \\

\begin{proposition}
Any $p \in \mathcal{P}(2)$ endows 
$Ho \spec_{BPU}$ with the structure of a symmetric monoidal category.
\end{proposition}
\textbf{Proof:}
First, recall that $BPU$ is a $K(\mathbb{Z},3)$-space.
We will make repeated use of this fact throughout the proof.\\
In Proposition \ref{symmstruct}, we 
have shown how to construct a symmetric monoidal structure on $Ho \specb$
from an action of a pointed $E_\infty$-operad on $B$. 
$\mathcal{P}$ is an $E_\infty$-operad and it acts on
$BPU$, but it is not pointed and 
therefore we can't apply Proposition \ref{symmstruct} directly.
However, the only parts of the proof where the pointedness of the operad
were used were
the construction of the unit isomorphism 
and the coherence diagramm involving it, so we will just reproof that.

We define
the bifunctor associated to $p \in \mathcal{P}$ to be
$$\underline{p}:=L (\theta(p)_! \circ \overline{\wedge}): 
Ho\spec_{BPU} \times Ho\spec_{BPU} \rightarrow Ho\spec_{BPU}$$
and the unit to be $\onebpu := i_{\ast !}S$.

To get the unit isomorphism, 
we need to find a homotopy from 
$$\xymatrix{ p_e: BPU \ar[rr]^{i_{\ast} \times id} & & BPU \times BPU
\ar[rr]^{\theta(p)} & & BPU}$$
 to the identity.
Note that 
$$\theta(p): BPU \times BPU \rightarrow BPU$$ induces the sum 
on $H^3$, i.e. $$\theta(p)_\ast = +: 
H^3(X,\mathbb{Z}) \times H^3(X,\mathbb{Z}) \rightarrow 
H^3(X,\mathbb{Z}) \ \ \ \forall X \in \Top$$
and thus  
$p_e$ is indeed homotopic to the identity.
Choosing such a homotopy defines by Proposition \ref{htpyisosb} 
a natural transformation
$\phi: L p_{e!} \stackrel{\sim}{\rightarrow} L id_!$.
This yields the required 
natural isomorphism 
$$\pp(\onebpu,X) \cong L p_{e!} X \cong QX \cong X.$$ Moreover,
the identity component of $\Top(BPU,BPU)$ is simply connected
(since $H^2(BPU,\mathbb{Z})=0$), so all homotopies connecting
$p_e$ and $id_{BPU}$ are homotopic and the isomorphism does not 
depend on the choice of the homotopy. 
It remains to check commutativity of
$$ 
\xymatrix{
\pp(\onebpu,\pp(X,Y)) \ar[r] \ar[d] & \pp(X,Y) \\
\pp(\pp(\onebpu,X),Y) \ar[ur]}.
$$
Recall that the associativity transformation was defined using
a path in the operad connecting
$\gamma(p;1,p)$ and $\gamma(p;p,1)$.
We denote the natural isomorphisms induced
by the various homotopies by $\phi_i$.
The commutativity
of the above diagram is the same as the commutativity of
$$\xymatrix{
\theta(p)_!(i_\ast \times id)_! \theta(p)_! (QX \ewedge QY) \ar[r]^/.5cm/{\sim} 
\ar[d]^{\sim}&
p_{e !} \theta(p)_! (QX \ewedge QY)  \ar[dd]^{\phi_0}  \\
\theta(p)_! (1 \times \theta(p))_! (i_\ast  \times id)_! (QX \ewedge QY) \ar[d]^{\phi_1}  \\
\theta(p)_! (\theta(p) \times 1)_! (i_\ast  \times id)_! 
(QX \ewedge QY) \ar[d]^\sim &  \theta(p)_! (QX \ewedge QY)\\
\theta(p)_! (p_e \times id)_! (QX \ewedge QY) \ar[ru]^{\phi_2}
 }$$
for all $X,Y \in \spec_{BPU}$. Note that 
$H^2(BPU \times BPU, \mathbb{Z}) = 0 $ and hence all homotopies between maps \\
$f,g: BPU \times BPU \rightarrow BPU$ are homotopic, so the application
of  Propositions \ref{htpyisosb}, \ref{funcsb} finishes the proof.
\hfill $\square$ 

\begin{proposition}
The gauge group $\mathcal{G}(EPU)$ of bundle automorphisms of the 
principal $PU$-bundle $EPU \rightarrow BPU$ is connected. 
\end{proposition}
\textbf{Proof:}
A gauge transformation $g \in \mathcal{G}(EPU)$ defines a 
principal $PU$-bundle $P_g \rightarrow BPU \times S^1$ by 
identifying $EPU \times \{0\}$ and $EPU \times \{1\}$ via $g$
in $EPU \times I$. Now, 
$H^3(BPU \times S^1, \mathbb{Z}) \cong H^3(BPU,\mathbb{Z})$ and therefore
$P_g$ is isomorphic to $\pi^{\ast} EPU$, where 
$$\pi:BPU \times S^1 \rightarrow BPU$$ is the projection. This implies that one can extend
$P_g$ to a $PU$-bundle on $BPU \times D^2$. The choice of such an extension
yields a path connecting $g$ and $id_{EPU}$ in $\mathcal{G}(EPU)$. 
  
\hfill $\square$ \\

\begin{lemma}
$E:=BPUK:=EPU \times_{PU} K$ is a monoid in
$(Ho \spec_{BPU},\underline{p})$ for all $p \in \mathcal{P}(2)$.

\end{lemma}
\textbf{Proof:} 
The non-pointedness of $\mathcal{P}$ obstructs the direct 
application of Lemma \ref{Monoid}, thus we will redo those parts of the proof where the basepoints of the operad were used, namely the construction of the unit
$\eta: \onebpu \rightarrow E$. 
Note that the map 
$ \tilde{\eta} : S \stackrel{\eta'} \rightarrow K' \stackrel{\phi}{\rightarrow} K$ 
is a homotopy unit for each $p \in \mathcal{P}(2)$ since
$$\xymatrix{
K'(V) \wedge K'(W) \ar[rr]^{\mu^+} \ar[d]^{\phi \wedge \phi}& 
& K'(V \oplus W)\ar[d]^{\phi} \\
K(V) \wedge K(W) \ar[r]^{\mu^+} & 
R^{+}_{H'(V\oplus W) \otimes \mathbf{H} \otimes \mathbf{H}} \ar[r] & 
K(V \oplus W)
}$$ is homotopy commutative.
Define $\eta:\onebpu = i_{\ast !} S \rightarrow E$ to be the adjoint of
$\tilde{\eta}$. The last thing to check is the commutativity
of
$$\xymatrix{
\pp(\onebpu,E) \ar[r] \ar[dr] & \pp(E,E) \ar[d] \\  
& E.}$$
By the definition of the unit isomorphism, all we have to find are maps 
$\psi_1, \psi_2$ that make the following diagram homotopy commutative.
In the proof of Lemma \ref{Monoid}, these were constructed from the action 
of $\mathcal{P}$. Since $p_e$ need not be of the form $\theta(q)$ for some 
$q \in \mathcal{P}(1)$, we cannot copy that construction. \\
\begin{sideways}
$$\xymatrix{
\theta(p)_!(Q1_{BPU} \ewedge QE) \ar[dr] \ar[d]^{\sim}
\ar[rr]^{\theta(p)_!(Q\eta \ewedge id)}
& & \theta(p)_!(QE \ewedge QE) \ar[dd]& \\
\theta(p)_!(1_{BG} \ewedge QE) \ar[r] \ar[d]^{\sim}
& \theta(p)_!(1_{BG} \ewedge E) \ar[dr]^{\theta(p)_!(\eta \ewedge id)} 
\ar[d]^{\sim} \ar@{}[ddddrr]^{\mathbf{I}}& & \\
\theta(p)_! (i_{\ast} \times id)_! (S \ewedge QE) \ar[d]^{\sim} \ar[r] &
\theta(p)_! (i_{\ast} \times id)_! (S \ewedge E) \ar[d]^{\sim} & 
\theta(p)_!(E \ewedge E) \ar[dddr]^{\mu(p)} & \\
p_{e !} QE \ar[r] \ar[d]^{\sim}& 
p_{e !} E \ar[ddrr]^{\psi_1} \ar[d]^{\sim} & & \\
h_! i_{0!}QE \ar[r] \ar[d]^{\sim} &
h_! i_{0!}E \ar[rd] \ar[d] & & \\
h_! \pi^{\ast} RQE \ar[r] & 
h_! \pi^{\ast} RE & 
h_! \pi^{\ast} E \ar[l]^{\sim} \ar[r]^{\psi_2} & E \\
h_! i_{1!}QE \ar[r] \ar[u]^{\sim} &
h_! i_{1!}E \ar[ru] \ar[u] & & \\
QE \ar[u]^{\sim} \ar[r]^{\sim} & E \ar[rruu]^{id} \ar[u]^{\sim}
}$$
\end{sideways}\\
To obtain $\psi_1$, let us first choose an isomorphism 
$\tilde{\psi_1}': EPU \rightarrow p_{e}^{\ast}EPU$ and use
this to define  $\tilde{\psi_1}: E \rightarrow p_{e}^{\ast} E.$
 

Letting 
$\psi_1 : p_{e !} E \rightarrow  E$ be the adjoint, 
the quadrangle $\mathbf{I}$ in the diagram commutes in $Ho\spec_{BPU}$ since 
$\tilde{\eta}$ is a homotopy unit.
Moreover, since the gauge group of $EPU$ is connected,
one can extend the isomorphism 
$\tilde{\psi_1}'$ to \\
$\tilde{\psi_2}': \pi^{\ast} EPU \rightarrow h^{\ast}EPU$ in such a way
that $i_{0}^{\ast}(\tilde{\psi_2}') = \tilde{\psi_1}'$ and 
$i_{1}^{\ast}(\tilde{\psi_2}') = id_{EPU}.$
Now, $\tilde{\psi_{2}}'$ yields 
$\tilde{\psi_2}: \pi^{\ast} E \rightarrow h^{\ast} E$. Defining $\psi_2$
to be the adjoint of $\tilde{\psi_2}$ completes the argument.  \hfill $\square$ \\
\section{Twisted $Spin^c$-cobordism}
To define twisted $Spin^c$-cobordism, we build upon the 
construction of  $M\spinc$ in \cite{Jo1}. Recall from the 
previous section the representation 
$$\varrho(V): Pin^{c}(V) \rightarrow \mathbf{U}_{H(V)}.$$ 
Restricting
it to the subgroup of even elements $Spin^{c}(V) \subset Pin^c (V)$, 
we obtain a representation $\spinc (V) \rightarrow \mathbf{U}^{0}_{H(V)}$ which we call
again $\varrho(V)$. 
Since $\mathbf{U}^{0}_{H(V)}$ is contractible, 
$\mathbf{U}^{0}_{H(V)} \times_{\spinc (V)} V$ is a model for the universal
$\spinc (V)$-bundle. We define the spaces of the spectrum as
$$M\spinc(V):= Th(\mathbf{U}^{0}_{H(V)} \times_{\spinc (V)} V) = 
\mathbf{U}^{0}_{H(V) +} \wedge_{\spinc (V)} S^V.$$ The
$O(V)$-action is obtained by combining the $O(V)$-action induced by
the conjugation action of $\spinc (V)$ on 
$\mathbf{U}^{0}_{H(V)}$ and the standard action on 
$S^V$. To define the structure maps
$$\sigma_{V,W}: M \spinc (V) \wedge S^W \rightarrow M \spinc (V \oplus W),$$ 
we use  the canonical isomorphism $H(V \oplus W) \cong H(V) \otimes H'(W)$
and the homeomorphism $f:S^V \wedge S^W \rightarrow S^{V \oplus W}$
and define $$\sigma_{V,W}[U,v,w]:=[U \otimes 1, f(v,w)].$$
Note that $PU=PU^0(\mathbf{H})$ acts on $M \spinc$ via
$$ PU \times M \spinc (V) \rightarrow M \spinc (V), \ \ \ \ \  \ \
([u],[U,v]) \mapsto [(1 \otimes u) U ,v].$$
This descends indeed to $PU$ as the central $U(1)$ of $Spin^c(V)$ acts trivially
on $V$ but as multiples of the identity on $U^0_{H(V)}$.
\begin{lemma}
$M\spinc$ is in $\mathcal{P}(\spec)$ and the
$PU$-action on $M\spinc$ is compatible with the $\mathcal{P}$-algebra
structures on $PU$ and $M\spinc$. 
\end{lemma}
 
Finally, there is an obvious map 
$$\alpha: M\spinc \rightarrow K, \
\alpha(U,v):= U \tilde{\eta}(v) U^{-1} \ \text{for} \ 
U \in \mathbf{U}^{0}_{H(V)}, v \in V$$
which is a $PU$-equivariant map of $\mathcal{P}$-spectra. 
\begin{proposition}\cite{Jo1}
$M\spinc$ represents $\spinc$-cobordism and $\alpha$ induces
the classical orientation homomorphism
$M Spin^{c \ast} \rightarrow K^\ast$.
\end{proposition}
Analogous to the $K$-spectrum over $BPU$, we define
an ex-spectrum over $BPU$ with fibre $M \spinc$ by
$\mathcal{M}_{\spinc} (V):= EPU \times_{PU} M\spinc$.
A similar proof as for the parametrised $K$-spectrum shows that
$\mathcal{M}_{\spinc}$ is a parametrised ring spectrum over $BPU$. 
\begin{definition}
The parametrised $M\spinc$-orientation is
$$\beta:=1 \times \alpha: EPU \times_{PU} M \spinc 
\rightarrow EPU \times_{PU} K.$$
\end{definition}


\begin{thebibliography}{99}
\bibitem{Adams} J.F. Adams,
\textsf{Stable Homotopy and Generalised Homology}, Chicago Lectures in Mathematics,
The University of Chicago Press, Chicago (1974)
\bibitem{AS1}
M. Atiyah and G.Segal,
\textsf{Twisted K-theory},
arXiv:math.KT/0407054 (2004)
\bibitem{AS2}
M. Atiyah and G.Segal,
\textsf{Twisted K-theory and cohomology},
arXiv:math.KT/0510674 (2005)
\bibitem{BJ}S. Baaj and P. Julg, \textsf{Th\'eorie bivariante de Kasparov et
op\'erateurs non born\'es dans les $C^\ast$-module hilbertiens}, C.R.
Acad. Sci. Paris S\'er. I Math.296:21 (1983)
\bibitem{Booth1}
P.I. Booth, \textsf{The exponential law of maps I}, Proc. London Math. Soc. (3) 20
(1970)
\bibitem{Booth2} 
P.I. Booth,
\textsf{The section problem and the lifting problem}
Math. Z. 121 (1971)
\bibitem{Booth3}
P.I. Booth,
\textsf{The exponential law of maps II},Math. Z. 121 (1971)
\bibitem{BrBooth1}
P.I. Booth and R. Brown, 
\textsf{Spaces of partial maps, fibred mapping spaces and the 
compact-open topology.}, General Topology and Appl. 8 (1978)
\bibitem{BrBooth2}
P.I. Booth and R. Brown, 
\textsf{On the application of fibred mapping spaces to exponential laws for bundles, ex-spaces and other categories of maps.}, General Topology and Appl. 8 (1978)
\bibitem{Brown}
E.H. Brown, \textsf{Abstract Homotopy Theory},
Trans. Amer. Math. Soc. 119 (1965)
\bibitem{BJS}
U.Bunke, M.Joachim and S.Stolz, \textsf{Classifying spaces and spectra representing the $K$-theory of a graded $C^ *$-algebra},
High-dimensional manifold topology,  80--102 (2003)

\bibitem{CaWang}A. L. Carey and B.-L. Wang
\textsf{Thom isomorphism and Push-forward map in twisted K-theory},
 arXiv:math.KT/0507414 (2005)


\bibitem{Clapp}
M. Clapp,
\textsf{Duality and transfer for parametrized spectra},
 Arch. Math. (Basel)  37, no. 5 (1981)
\bibitem{ClappPuppe}
M. Clapp and D. Puppe,
\textsf{The homotopy category of parametrized spectra},  
Manuscripta Math.  45,  no. 3  (1984)

\bibitem{DS} W.G Dwyer and J. Spalinski,
\textsf{Homotopy theories and model categories}
Handbook of algebraic topology, North-Holland, Amsterdam (1995) 
\bibitem{Dyer} E.Dyer, \textsf{Cohomology Theories} Mathematical Lecture Notes Series,
W.A. Benjamin Inc, New York, (1969)
\bibitem{FHT1}D. S. Freed, M. J. Hopkins and C.
        Teleman,
\textsf{Twisted equivariant K-theory with complex coefficients},
arXiv:math.AT/0206257 (2002)
\bibitem{FHT2}D. S. Freed, M. J. Hopkins and C.
        Teleman,
\textsf{Twisted K-theory and loop group representations},
arXiv:math.AT/0312155 (2003)
\bibitem{FHT3}D. S. Freed, M. J. Hopkins and C.
        Teleman,
\textsf{Loop Groups and Twisted K-Theory II},
arXiv:math.AT/0511232 (2005)


\bibitem{Hovey} M. Hovey,
\textsf{Model Categories.} Amer. Math. Soc. Surveys and Monographs 63 (1999) 
\bibitem{Intermont}
M. Intermont, M.W. Johnson,
\textsf{ Model structures on the category of ex-spaces}
Topology Appl. 119 (2002)


\bibitem{James1}
I.M. James,
\textsf{Fibrewise topology},
Cambridge Tracts in mathematics 91, Cambridge University Press, Cambridge (1989)
\bibitem{Jo1}
M. Joachim, \textsf{Higher coherences for equivariant K-Theory}
Structured Ring Spectra,
London Math. Soc. Lecture Notes Ser. 315, Cambridge University Press, 
Cambridge (2004)

\bibitem{Kelly}
G.M.Kelly,
\textsf{On MacLane's Conditions for Coherence of Natural Associativities,
Commutativities, etc.}
Journal of Algebra 1 (1964)
\bibitem{MM}
M.A. Mandell and J.P. May,
\textsf{Equivariant orthogonal spectra and $S$-modules},
Mem. Amer. Math. Soc.  159, no. 755 (2002)

\bibitem{MMSS}
M.A. Mandell, J.P. May, S. Schwede and B. Shipley,
\textsf{Model categories of diagram spectra},
Proc. London Math. Soc (3) 82 (2001)

\bibitem{Markl}
M. Markl, S. Shnider and J. Stasheff, 
\textsf{Operads in Algebra, Topology and Physics}
Mathematical Surveys and Monographs, 96,
AMS, Providence (2002) 
\bibitem {GeomILS}
J.P. May \textsf{The Geometry of Iterated Loop Spaces}, 
Lecture Notes in Mathematics, Vol. 271, Springer-Verlag, Berlin-New York (1972)
\bibitem{Mayarticle}
J.P. May \textsf{$E_\infty$ spaces, group completions, and permutative categories}
New developments in topology (Proc. Sympos. Algebraic Topology, Oxford, 1972),  London Math. Soc. Lecture Note Ser., No. 11, Cambridge Univ. Press, London (1974)
\bibitem{MAddi}
J.P. May, \textsf{The additivity of traces in triangulated categories}
 Adv. Math. 163,  no. 1, (2001)

\bibitem{MS}
J.P. May and J. Sigurdsson,
\textsf{Parametrized Homotopy Theory},
arXiv:math.AT/0411656 (2004)
\bibitem{Schoen}
R. Sch\"on,
\textsf{ On representable half-exact functors over $B$,}
Quaestiones Math. 13, no. 3-4, (1990)

\bibitem{Sw}
R.M. Switzer, \textsf{Algebraic Topology - Homotopy and Homology},
Die Grundlehren der mathematischen Wissenschaften 212, Springer-Verlag 
New York-Heidelberg (1975)

\bibitem{TXL1}J.-L. Tu, P. Xu and C. Laurent-Gengoux
\textsf{Twisted K-theory of differentiable stacks},
arXiv:math.KT/0306138 (2003)

\bibitem{TX}J.-L. Tu and P. Xu
\textsf{The ring structure for equivariant twisted K-theory},
arXiv:math.KT/0604160 (2006)

\bibitem{Quillen}
D.G Quillen, \textsf{Homotopical algebra}, lecture notes in Mathematics 43,
Springer-Verlag New York-Heidelberg (1967)
\end{thebibliography}
\end{document}